\RequirePackage{fix-cm}
\documentclass{svjour3}                     
\smartqed  
\usepackage{graphicx}
\usepackage[top=1.1in, bottom=1.1in, left=1.2in, right=1.2in]{geometry}
\usepackage{lineno,hyperref}
\usepackage{bbm}
\usepackage{amssymb}
\usepackage[misc]{ifsym}
\usepackage{multirow,booktabs}
\usepackage{caption}
\usepackage{mathrsfs}
\usepackage{amsfonts}
\usepackage{colortbl,dcolumn}
\usepackage{tabularx}
\usepackage[tbtags]{amsmath}
\usepackage{verbatim}
\usepackage{xspace}
\usepackage{booktabs}
\usepackage{color}
\newcommand{\R}{\mathbb{R}}

\newcommand{\E}{\mathbb{E}}
\renewcommand{\P}{\mathbb{P}}
\newcommand{\dd}{\text{d}}
\newcommand{\Ftau}{F^{(\tau)}}

\newcommand{\RefBEM}{\texttt{RefBEM}\xspace}
\newcommand{\Splitting}{\texttt{Splitting}\xspace}
\newcommand{\BEM}{\texttt{BEM}\xspace}
\newcommand{\TSM}{\texttt{TSM}\xspace}
\newcommand{\TEM}{\texttt{TEM}\xspace}
\usepackage[numbers,sort&compress]{natbib}
\begin{document}

\title{Convergence of a splitting method for a general interest rate model }


\titlerunning{Tamed-Spiltting scheme for Ait-Salahia model}        

\author{ Gabriel Lord \and  Mengchao Wang \Letter    }
\authorrunning{ Gabriel Lord  \and
	  Mengchao Wang  } 
\institute{ 
	Gabriel Lord   \at
	Department of Mathematics, IMAPP, Radboud University, Nijmegen, The Netherlands\\
	\email{gabriel.lord@ru.nl}
     \and
     Mengchao Wang 
     \Letter  \at
	Department of Mathematics, IMAPP, Radboud University, Nijmegen, The Netherlands\\
	School of Mathematics and Statistics, HNP-LAMA, Central South University, Changsha 410083, P. R. China \\
	\email{mengchao.wang@ru.nl}
	\and
	This work was supported by China Scholarship Council,
 Graduate Research and Innovation Program of Central South University(No.CX20220246)
	and Radboud University Nijmegen. 
}
\date{}
\maketitle
\begin{abstract}
  We prove mean-square convergence of a novel numerical method, the tamed-splitting method, for a generalized Ait-Sahalia interest rate model. The method is based on a Lamperti transform, splitting and applying a tamed numerical method for the nonlinearity.
The main difficulty in the analysis is caused by
	the non-globally Lipschitz drift coefficients of the model.
        We examine the existence, uniqueness of the solution and boundedness of moments for the transformed SDE. We then prove bounded moments and inverses moments for the numerical approximation.
        The tamed-splitting method is a hybrid method in the sense that a backstop method is invoked to prevent solutions from overshooting
zero and becoming negative.
We successfully recover the mean-square convergence rate of order one for the tamed-splitting method. 
In addition we prove that the probability of ever needing the backstop method to prevent a negative value can be made arbitrarily small.
In our numerical experiments 
we compare to other numerical methods in the literature for realistic parameter values.

\end{abstract}
\section{Introduction} \label{sect:introduction}
Stochastic differential equations (SDEs) are widely used in various scientific areas to 
model real-life phenomena affected by random noise.
Here we consider the generalized Ait-Sahalia interest rate model introduced in \cite{ait1996testing} 
\begin{equation}\label{eq:Ait-Sahalia-model}
	\dd X_t = (a_{-1} X_t^{-1} - a_{0}
	+ a_{1} X_t - a_{2} X_t^{\gamma})
	\, \dd t
	+b X_t^{\theta} \dd W_t , \quad
	X_0=x_0, 
	\quad
	t > 0
\end{equation}
which has since been investigated  by various authors (see, e.g. \cite{conley1997short}, \cite{gallant1997estimation}, \cite{hong2005nonparametric}). 
Here, $a_{-1}$,\,$a_{0}$,\,$a_{1}$,\,$a_{2}$,\,$ b > 0$, both $\theta$,\,$\gamma>1$ and $ \{W_t\}_{t \in [0,\infty)}$ is a one-dimensional Brownian motion on a filtered probability space
$ (\Omega, \mathscr{ F }, \P, \{\mathscr{ F }_t\}_{t \geq 0}) $
with respect to the normal filtration $  \{\mathscr{ F }_t\}_{t \geq 0} $.

The model \eqref{eq:Ait-Sahalia-model} 
clearly violates the Lipschitz and linear growth conditions which are traditionally imposed in the study of SDEs and their simulation (see \cite{kloeden_numerical_2011,lord2014introduction,giles2008multilevel}). As such it forms an interesting test case for numerical methods, 
for $\theta\in[\frac12,1]$ see for example \cite{alfonsi2005discretization,berkaoui2008euler,bossy2007efficient,higham2005convergence}.
For the diffusion coefficient $ \theta > 1 $, as we consider here,
the model was studied in \cite{Szpruch2011Numerical},
where strong convergence of the backward Euler method was proved,
but without revealing any convergence rate.  Deng et.al. \cite{Deng2018Generalized} examined
the analytical properties of the model with Poisson jumps, including positivity,
boundedness and pathwise asymptotic estimates. They applied the Euler-Maruyama (EM) method 
and proved that the explicit scheme converges in probability to the true solution of the model.
Emmanuel et.al in \cite{emmanuel2021truncated} study
analytical properties for the true solution of 
the model with delay and construct new truncated Euler-Maruyama (EM) methods to study properties of the numerical solutions under the local Lipschitz condition plus a Khasminskii-type condition.
Zhao et.al. in  \cite{zhao2020backward} proved a 
mean-square convergence rate of order one half for the backward Euler method (BEM) for a generalized Ait-Sahalia interest rate 
model with Poisson jumps and showed that the BEM preserves 
the positivity of the original problem.

To the best of our knowledge, no strong convergence order 1 has been reported 
in the literature for explicit numerical approximations of \eqref{eq:Ait-Sahalia-model}.
In this paper, we show that a strongly convergent numerical scheme
can be constructed by an application of the Lamperti transform followed
by a splitting approach where, rather than the nonlinearities being solved exactly, a tamed method such as in \cite{hutzenthaler2012} is used. As such our analysis combines splitting approaches with a taming scheme. We call this the tamed-splitting method (\TSM).
In our numerical experiments we compare the tamed-splitting method to a standard splitting method (solving the drift nonlinearities exactly, denoted \Splitting), the backward Euler method (\BEM) for the Lamperti transform, tamed Euler
method (\TEM) for the transformed SDE, and backward Euler method for the original SDE \eqref{eq:Ait-Sahalia-model} (\RefBEM).

%
The structure of the article is as follows. 
In Section \ref{sect:Ait-Sahalia} we present the SDE arising from the Lamperti transform of \eqref{eq:Ait-Sahalia-model} and examine  existence, uniqueness as well as prove moment and inverse moment bounds on the solution. 
In Section \ref{sec:derivation}, we introduce our two numerical methods. For the tamed-splitting method we prove moment and inverse moment bounds of the numerical solution in Section \ref{sec:convergence}. Then, by introducing an auxiliary process, we prove the mean squared convergence rate of the 
tamed-splitting method is one.
In Section \ref{sec:positivity} we show that the probability of a numerical solution taking a negative value is arbitrarily small,
that is, this probability tends to zero as the step size tends to zero.
Finally in Section \ref{sec:numerics} we numerically compare
convergence and efficiency of several commonly used methods.

\section{The Ait-Sahalia model} \label{sect:Ait-Sahalia}
We introduce the notation we use throughout.
Given a filtered probability space
$ ( \Omega, \mathcal{ F }, \{\mathcal{ F }_t\}_{t \in [0,T]} ,\P ) $
satisfying the usual hypotheses we let $W_t$ be a Brownian motion defined on that space.
We let $ \E $ denote the expectation and $L^p(\Omega;\R) $ the space of 
$p$-times integrable random variables with 
$ \| \xi \|_{L^p(\Omega;\R)} := \big ( \E [ |\xi |^p ] \big)^{ 1/p}$ for any $ p \geq 1 $.
Let 
$ x \wedge y :=\min \{x, y\} $
for any $ x, y \in \mathbb{R} $. For notational simplicity,
the letter $ C $ is used to
denote a generic positive constant, which is independent of the
time stepsize and may vary for each appearance.
%
%
%
%

The well-posedness of the  Ait-Salalia model \eqref{eq:Ait-Sahalia-model} was
proved in \cite[Theorem 2.1]{Szpruch2011Numerical}.
\begin{proposition} \label{prop:Ait-Sahalia-solution-well-posedness}
	Let the initial data $ X_0 = x_0 >0 $.
	For constants $ a_{-1},\, a_{0}, \, a_{1}, \, a_{2}, \, b >0 $ 
	and $ \gamma, \, \theta > 1 $,
	the problem \eqref{eq:Ait-Sahalia-model}
	admits a unique positive global solution.
\end{proposition}
We now apply the Lamperti transform to obtain an SDE with additive noise. Taking  $Y = X^{1-\theta} $ and using It\^o's formula we get %
\begin{equation}\label{eq:lamperti-translation1-jump}
d Y_t	= f(Y_t) \dd t + b(1-\theta) \dd W_t,  \quad  t \in(0, T],
\end{equation}
where
\begin{equation}\label{eq:defnf}
f(x) := (\theta-1) 
 \big( a_2 x^{\frac{\gamma-\theta}{1-\theta}}
 - a_1 x - 
a_{-1} x^{\frac{\theta+1} {\theta-1}} 
+  a_0 x^{\frac{\theta}{\theta-1}}
+ \frac{1}{2} b^2 \theta  x^{-1} \big)
\end{equation}
and $  Y_0=X_{0}^{1-\theta} \in \mathbb{R}^{+}$. We first examine existence and uniqueness for the transformed SDE \eqref{eq:lamperti-translation1-jump}. 
Before proceeding furthermore,
we introduce the operator
$\mathbb{L} : C^{2} \big( \R_{+} \times [ 0, \infty ) ,\R \big) \rightarrow 
C \big( \R_{+} \times [ 0, \infty ) ,\R \big)$ 
defined by
\begin{equation} \label{eq:diffusion-operator}
	\mathbb{L} \phi (x,t) := \phi_{x}(x,t) 
	f(x) + \tfrac{1}{2} \phi_{xx}(x,t) (b (1 - \theta ) )^{2},
	\quad
	\phi (x,t) \in C^{2} \big( \R_{+} \times [ 0, \infty ) ,\R \big).
\end{equation}
If the second variable $ t $ vanishes, 
we can rewrite \eqref{eq:diffusion-operator}
as follows
\begin{equation*}
	\mathbb{L} \phi (x) := \phi_{x}(x) 
	f(x) + \tfrac{1}{2} \phi_{xx}(x) (b (1 - \theta ) )^{2},
	\quad
	\phi (x) \in C^{2} \big( \R_{+} ,\R \big).
\end{equation*}
\begin{proposition} \label{prop:Lamperti-solution-well-posedness}
Let the initial data 
	$ Y_0=X_0^{1-\theta} = x_0^{1-\theta} >0 $.
	For constants $ a_{-1},\, a_{0}, \, a_{1}, \, a_{2}, \, b >0 $ 
	and $ \gamma, \, \theta > 1 $,
	the problem 
	\eqref{eq:lamperti-translation1-jump}
	admits a unique positive global solution, which 
	almost surely satisfies
	\begin{equation} \label{eq:solution-Lamperti-Ait-Sahalia-mode}
		Y_t = Y_0 + \int_{0}^{t}
		f( Y_s )  \, \dd s
		+\int_{0}^{t} b(1-\theta)  \dd  W_s,
		\quad
		t \geq 0.
	\end{equation}
\end{proposition}
\textbf{Proof of Proposition 
	\ref{prop:Lamperti-solution-well-posedness}}
It is straightforward to see 
that the drift 
coefficients of 
\eqref{eq:lamperti-translation1-jump} 
are locally Lipschitz continuous in 
$(0, \infty)$. Following the standard arguments in
\cite{Mao2008Stochastic}  and noting $ Y_0 > 0 $,
one can show that there is a  unique maximal local solution 
$Y_t, \, t \in\left[0, \tau_{e}\right)$, 
where $\tau_{e}$ is the stopping time of 
the explosion or first zero time.
To confirm we have a global solution, 
we need to prove $ \tau_{e} = \infty $ a.s.
For any sufficiently large positive integer $n$, 
satisfying $1/n<Y_0<n$, we define the stopping times
\begin{equation}\label{eq:defn-stopping-times}
	\tau_{n} := \inf \left\{t \in\left[0, 
	\tau_{e}\right): Y_t \notin(1 / n, n)\right\},
\end{equation}
where throughout this paper we set 
$ \inf ( \emptyset ) = \infty $.
Obviously $ \tau_n $ is increasing as $ n \rightarrow \infty $
and we set $\tau_{\infty} := \lim_{n \rightarrow \infty} \tau_{n}$. 
In view of \eqref{eq:defn-stopping-times}, one knows $ \tau_{\infty} \leq \tau_{e} $ a.s.
If we can prove $\tau_{n} \rightarrow \infty$ a.s.
as $n \rightarrow \infty$, 
then $\tau_{e}=\infty$ a.s. 
and $Y_t > 0$ a.s. 
for all $t \geq 0$, then the proof is complete.
To prove $\tau_{\infty}=\infty$ a.s., it suffices to show that 
$\P\left\{\tau_{n} \leq T\right\} \rightarrow 0$ 
as $n \rightarrow \infty$ 
for any constant $T>0,$ which immediately implies $ \P \left\{\tau_{\infty}=\infty\right\}=1$, as required. \\
Given a fixed constant  $ \alpha \in (0,1)$,
let us define a function $V \in C^{2}\big((0,\infty), (0,\infty)\big)$ by
\begin{equation} \label{eq:defn-V}
	V(x):=x^{\alpha}-\alpha \log x.
\end{equation}
It is easy to check that $V(x) \rightarrow \infty$ 
as $ x \rightarrow \infty$ or $x \rightarrow 0 $ and that
\begin{equation}
	\label{eq:VxVxx-computation}
	V_{x}(x)=\alpha\big(x^{\alpha-1}-x^{-1}\big),
	\quad
	V_{xx}(x)=\alpha(\alpha-1) x^{\alpha-2}+\alpha x^{-2}.
\end{equation}
%
Bearing \eqref{eq:VxVxx-computation} in mind, we see that 
	\begin{align}  
	\label{eq:LV-computation}
		 \mathbb{L} V(x)   
		&  =
		\alpha \left(x^{\alpha-1}-x^{-1}\right)
		( \theta - 1 )
		\left[ a_2 x^{ \frac{\gamma-\theta}{1-\theta} }
		- a_1 x - a_{-1} x^{\frac{\theta+1} {\theta-1}} 
		+  a_0 x^{\frac{\theta}{\theta-1}}
		 + \frac{1}{2} b^2 \theta  x^{-1} \right] \nonumber \\
		& \quad 
		+ \frac{1}{2} b^2 (1 - \theta )^2 \left[\alpha(\alpha-1) 
		x^{\alpha-2}+\alpha x^{-2}\right] \nonumber \\
		 &  =
		 \alpha ( \theta - 1 )
		 \left( a_2 x^{ \alpha-1 + 
		 	\frac{\gamma-\theta}{ 1 - \theta } }
		 - a_1 x^{\alpha} - a_{-1} 
		 x^{\alpha-1+ \frac{\theta+1} {\theta-1}} 
		 +  a_0 x^{\alpha-1 + \frac{\theta}{\theta-1}}
		 +  \frac{1}{2} b^2 \theta  x^{\alpha - 2}
		  \right) \nonumber \\
		  & \quad 
		  - \alpha ( \theta - 1 )
		  \left( a_2 x^{-1 + 
		  	\frac{\gamma-\theta}{ 1 - \theta } }
		  - a_1  - a_{-1} 
		  x^{  \frac{ 2 } {\theta-1}} 
		  +  a_0 x^{ -1 + \frac{\theta}{\theta-1}}
		   + \frac{1}{2} b^2 \theta  x^{ - 2}
		  \right) \nonumber \\
		 & \quad 
		 + \frac{1}{2} b^2 (1 - \theta )^2 \left[\alpha(\alpha-1) 
		 x^{\alpha-2}+\alpha x^{-2}\right]  .  
	\end{align}
Taking \eqref{eq:LV-computation} 
into account and recalling 
$0<\alpha<1, \gamma>1$, $\theta>1$, $\lambda > 0$,
there is a constant $ K_1 > 0$ such that 
\begin{equation}\label{eq:bound-LV}
	\sup_{x \in(0, \infty)}
	\mathbb{L} V(x)
	\leq K_1
	< \infty.
\end{equation}
Indeed, if $ \gamma + 1 > 2 \theta $,
we can directly get 
$ \alpha-1 + \frac{\theta+1} {\theta-1}
 >  \frac{ 2 } {\theta-1} $. 
Thus, it is easy to see that the highest power of 
$ x $ is 
$ \alpha-1 + \frac{\theta+1} {\theta-1} $ 
and the lowest power of 
$ x $ is 
$ -1 +  \frac{\gamma-\theta}{ 1 - \theta } $ 
in \eqref{eq:LV-computation}. 
As a result of 
$ \alpha a_{2} ( \theta - 1 ) > 0 $
and
$ \alpha a_{-1} ( \theta - 1 ) > 0 $, 
there is a constant 
$ K_{1} > 0 $ such that 
\eqref{eq:bound-LV}
is fulfilled. 
If $ \gamma + 1 = 2 \theta $, 
it is easy to see that the highest power of 
$ x $ is 
$ \alpha - 1 + \frac{\theta+1} {\theta-1} $ 
and the lowest power of 
$ x $ is 
$ -2  $ 
in \eqref{eq:LV-computation}. 
As a result of 
$ \alpha a_{ -1 } ( \theta - 1 ) > 0 $
and
$ \frac{1}{2}  \alpha b^2 (1 - \theta )^2
- \alpha ( \theta - 1 ) a_{ 2 }
- \frac{1}{2}  \alpha b^2 \theta
 ( \theta - 1 ) < 0 $, 
there is a constant 
$ K_{1} > 0 $ such that 
\eqref{eq:bound-LV}
is fulfilled. 
By the It\^{o} formula \cite{Gardon2004approximations} 
applied to 
$ V ( Y_{ t \wedge \tau_{n} } ) $, 
$ t \in [ 0, T]$, 
we infer
\begin{equation*}
	\mathbb{E} 
	[
	V ( Y_{ T \wedge \tau_{n}} )
	]
	\leq V(Y_0 ) + K_1 T
	< \infty,
	\quad
	\text{for all} \, T >0.
\end{equation*}
Owing to the definitions \eqref{eq:defn-stopping-times}, \eqref{eq:defn-V}, we deduce from the above estimate that
\begin{equation}
	\mathbb{P}\left(\tau_{n} \leq T\right)[V(1/n) \wedge V(n)] 
	\leq \mathbb{E} V(Y_{ T \wedge \tau_{n}} ) 
	\leq V\left( Y_0 \right)+ K_1 T < \infty.
\end{equation}
This implies that 
$ \lim_{n \rightarrow \infty}
\mathbb{P}\left(\tau_{n} \leq T\right) = 0$ for any constant $T>0$ and the proof is thus complete.
\qed
In the following error analysis, the moment bounds of the solution to \eqref{eq:lamperti-translation1-jump} are frequently used.
We first prove when the $p$-th negative moments of the solution to \eqref{eq:lamperti-translation1-jump} are bounded.
\begin{lemma} \label{lem:solution-moment-bound}
	Let all conditions in Proposition \ref{prop:Ait-Sahalia-solution-well-posedness} 
	hold and let $ \{ Y_t \}_{ t \geq 0 } $ 
	be the unique solution to 
	\eqref{eq:lamperti-translation1-jump}, 
	given by \eqref{eq:solution-Lamperti-Ait-Sahalia-mode}.
	If one of the
	following two conditions holds: \\
	(i) $ p \geq \frac{ 2 }{ \theta - 1 } $ 
	when $ \gamma+1> 2\theta $; \\
	(ii) $ \frac{ 2 }{ \theta - 1 }
	 \leq p \leq  \frac{2 a_2 + b^2}{(\theta-1)b^2} $ 
	 when $ \gamma+1=2\theta $,\\
	then
	\begin{equation}\label{eq:p-moment-bound}
		\sup _{t \in[0, \infty)} \E [ |Y_t |^{-p} ]<\infty.
	\end{equation}
\end{lemma}
\textbf{Proof of Lemma \ref{lem:solution-moment-bound}}
For a sufficiently large positive integer
$ n $ satisfying $ \frac{1}{n}< x_0 < n $,
we define the stopping time
\begin{equation}\label{eq:tau_n}
  \tau_{n} := \inf  \{t \in [0, \infty ) :
    Y_t \notin(1 / n, n) \}.
\end{equation}
Also, we define
$ V_1 $ : $ \mathbb{R}_{+} 
\times [0,\infty) \rightarrow \mathbb{R}_{+} $
as follows
\begin{equation*}
	V_1(x,t) : = e^{t} x^{-p}, 
	\quad x \in \mathbb{R}_{+}, \, t \in [0,\infty).
\end{equation*}
We compute that
\begin{align}\label{eq:operator-L-computation}
		\mathbb{L} V_1(x,t) 
		& = e^{t} \big[ - p (\theta-1) x^{-p-1}
		(a_2 x^{\frac{\gamma-\theta}{1-\theta}}
		- a_1 x - 
		a_{-1} x^{\frac{\theta+1} {\theta-1}} 
		-  a_0 x^{\frac{\theta}{\theta-1}}
		+ \frac{1}{2} b^2 \theta x^{ - 1 } )  \nonumber \\
		& \quad
		+ \frac{1}{2} b^2 (1-\theta)^2 p (p+1)  
		x^{-p-2} \big]
		 \nonumber \\
			& = e^{t} \big[  p (\theta-1) 
		(-a_2 x^{-p-1+\frac{\gamma-\theta}{1-\theta}}
		+ a_1 x^{-p} +
		a_{-1} x^{-p-1+\frac{\theta+1} {\theta-1}} 
		+  a_0 x^{-p-1+\frac{\theta}{\theta-1}} \nonumber \\
			& \quad 
		+ \frac{1}{2} b^2 ( \theta - 1 )
		 p ( p ( \theta - 1 ) - 1 )  x^{-p-2} )
		\big] ,
\end{align}
where  the operator
$ \mathbb{L}  $
is defined in \eqref{eq:diffusion-operator}.
Infact,  the lowest power of 
$ x $ are 
$ - p - 1 +  \frac{\gamma-\theta}{ 1 - \theta } $ 
for 
$ \gamma + 1 > 2 \theta $
and
$ - p - 2 $
for
$ \gamma + 1 = 2 \theta $
in \eqref{eq:operator-L-computation},
respectively.
Also,
it is easy to see that the highest power of 
$ x $ is 
$ - p - 1 + \frac{\theta+1} {\theta-1} $ 
in \eqref{eq:operator-L-computation}.
Since $p \geq \frac{ 2 }{ \theta - 1 }$,
we can get
$ - p - 1 +  \frac{ \theta + 1 }{ \theta - 1 } \leq 0 $.
If $ \gamma + 1 = 2 \theta $, due to condition $(ii)$, we can obtain that
\begin{equation*}
  \frac{1}{2} b^2 ( \theta - 1 )
  p ( p ( \theta - 1 ) - 1 )
  - p ( \theta - 1 ) a_{ 2 } \leq 0.
\end{equation*}
Therefore, one can find a constant $ K_2 > 0$, such that
\begin{equation*}\label{ineq:operator-L-bound}
  \mathbb{L} V_1(Y_t)  \leq K_2 e^t.
\end{equation*}
By the It\^{o} formula for any $ t \geq 0 $,
\begin{equation*}
  \E[e^{t \wedge \tau_{n}} Y_{t \wedge \tau_{n}}^{-p}]
  \leq Y_{0}^{-p} + K_2 e^{t}.
\end{equation*}
Letting $n  \rightarrow \infty$ and applying Fatou's lemma,
we obtain
\begin{equation*}
  \E [|Y_t|^{-p}] \leq \frac{x_0^{-p}}{e^{t}} + K_2, \quad t \geq 0.
\end{equation*}
The proof of Lemma \ref{lem:solution-moment-bound} is thus completed.
\qed
Having obtained negative moment bounds the following lemma gives positive moment bounds of the solution to \eqref{eq:lamperti-translation1-jump}. 
\begin{lemma}\label{lem:solution-inverse-moment-bound}
  Let all conditions in Proposition \ref{prop:Ait-Sahalia-solution-well-posedness} hold with $ \gamma +1 \geq 2 \theta $ and let $ \{ Y_t \}_{ t \geq 0 } $  
  be the unique solution to \eqref{eq:lamperti-translation1-jump}, 
  given by \eqref{eq:solution-Lamperti-Ait-Sahalia-mode}.
  Then for any $ p \geq (\gamma  - 1 )/(\theta - 1 ) $
  it holds that
  \begin{equation*}\label{eq:inverse-p-moment-bound}
    \sup _{t \in[0, \infty)} \E [ |Y_t |^{p} ]<\infty.
  \end{equation*}
\end{lemma}
\textbf{Proof of Lemma \ref{lem:solution-inverse-moment-bound}.}
Define
$ V_2 $ : $ \mathbb{R}_{+} 
\times [0,\infty) \rightarrow \mathbb{R}_{+} $
as follows
\begin{equation}
  V_2 ( x, t) := e^{t} x^{p}, 
  \quad x \in \mathbb{R}_{+}, \,\,\, t \in [0,\infty).
\end{equation}
Here, $ \tau_{n}$ is defined in \eqref{eq:tau_n} and the operator $ \mathbb{L} $ in \eqref{eq:diffusion-operator}.
Then we have
\begin{align}\label{eq:operator-L-computationoperator-inverde-moment}
		\mathbb{L} V_2(x,t) 
		& = e^{t} \big[  p (\theta-1) x^{p-1}
		(a_2 x^{\frac{\gamma-\theta}{1-\theta}}
		- a_1 x - 
		a_{-1} x^{\frac{\theta+1} {\theta-1}} 
		+  a_0 x^{\frac{\theta}{\theta-1}}
		 +
		 \frac{ 1 }{ 2 } b^2 \theta x^{ - 1 }
		 )  \nonumber \\
		& \quad 
		+ \frac{1}{2} b^2 (1-\theta)^2 
		p (p - 1)  x^{p-2} \big]
		 \nonumber \\
		& = e^{t} \big[  p (\theta-1) 
		(a_2 x^{p-1+\frac{\gamma-\theta}{1-\theta}}
		- a_1 x^{p} -
		a_{-1} x^{p-1+\frac{\theta+1} {\theta-1}} 
		+  a_0 x^{p-1+\frac{\theta}{\theta-1}} \nonumber \\
		& \quad
		+ \frac{1}{2} b^2 ( \theta - 1 )
		 p  [ ( \theta - 1 )(p-1) + \theta ]  x^{p-2} )
		\big]  .
\end{align}
The highest power of
 $ x $
is
$ p-1+\frac{\theta+1} {\theta-1} $
in \eqref{eq:operator-L-computationoperator-inverde-moment}.
The coefficient of this term is negative.
According to 
$ p \geq \frac{ \gamma  - 1 }{  \theta - 1 }$,
the lowest power of $ x $
is non-negative, i.e.
$ p - 1 + \frac{\gamma-\theta}{1-\theta} \geq 0 $.
Therefore,
there exists a constant $ K_3 > 0 $ such that
\begin{equation}\label{ineq:operator-L-bound-inverse-moment}
	\mathbb{L} V_2(X_t)  \leq K_3 e^t.
\end{equation}
The remaining proof is similar to 
the proof of Lemma
\ref{lem:solution-moment-bound}
and thus omitted.
\qed

\section{Derivation of the splitting and tamed-splitting methods}
\label{sec:derivation}
In this section, we first introduce the splitting and then the tamed-splitting method for the 
Ait-Sahalia model \eqref{eq:Ait-Sahalia-model} which are based on the additive noise
SDE \eqref{eq:lamperti-translation1-jump}.
We then prove moment bounds and strong convergence rate of the  tamed-splitting method.

For  $ N \in \mathbb{N} $ we construct on $[0,T]$ a uniform mesh $t_{i+1}=t_i+\tau$, $t_0=0$
with $\tau =\frac{T}{N}$ being the stepsize and introduce the notation
\begin{equation*}
  \lfloor t \rfloor : = t_{i}, \quad 
  \text{for} \quad
  t \in [ t_{i} , t_{i+1} ), \quad 
    i \in \{ 0, 1, \cdots, N-1 \}.
\end{equation*}
We can re-write \eqref{eq:lamperti-translation1-jump} as
\begin{equation}\label{eq:Lamperti}
  dY_t = [-\lambda Y_t + F(Y_t)] dt + b(1-\theta) dW_t,
  \end{equation}
where $f(x) = -\lambda x + F(x)$ with 
$ \lambda := (\theta-1)a_1 $
and
\begin{equation}\label{eq:defLamF}
  F (x) 
:= (\theta-1)\big[ a_2 x^{\frac{\gamma-\theta}{1-\theta}}
- a_{-1} x^{\frac{\theta+1} {\theta-1}} 
+ a_0 x^{\frac{\theta}{\theta-1}} 
+ \frac{1}{2} b^2 \theta x^{-1} \big].
\end{equation}
By the mean value theorem of differentiation
\begin{equation}\label{eq:controal_F}
 \begin{split}
 \big( F(x) - F(y) \big) \cdot
 ( x - y ) 
 & =
 F'(\xi
 ) | x - y |^2  \\
 & =
 \big( - a_2(\gamma-\theta) \xi^{\frac{\gamma-1}{1-\theta}}
- a_{-1} (\theta+1)
\xi^{\frac{2} {\theta-1}} 
+ a_0 \theta
\xi^{\frac{1}{\theta-1}} 
- \frac{\theta-1}{2} b^2 \theta \xi^{-2} \big) 
| x - y |^2 \\
& \leq C | x - y |^2,
 \qquad
 \xi \in (x,y)
 \,\,\,
 \text{or}
 \,\,\,
 \xi \in (y,x),
 \quad  x,y \in \R_{+}. 
  \end{split}
\end{equation}
The (standard) splitting method is derived from \eqref{eq:Lamperti} by solving the four ODEs arising from the  nonlinearity $F$ along with the OU process. That is we solve
the systems
\begin{equation*}
\left\{\begin{array}{l}
	\mathrm{d} U_t=(\theta-1) a_{2} U_t^{\frac{\gamma-\theta}{1-\theta}} \mathrm{d} t \\
	\mathrm{~d} V_t=(1- \theta ) a_{-1} V_t^{\frac{\theta+1}{\theta-1}} \mathrm{~d} t \\
	\mathrm{~d} P_t=(\theta - 1) a_{0} P_t^{\frac{\theta}{\theta-1}} \mathrm{~d} t \\
	\mathrm{~d} Q_t=\frac{1}{2} b^{2} \theta(\theta-1) Q_t^{-1} \mathrm{~d} t \\
	\mathrm{~d} Z_t= \lambda Z_t \mathrm{~d} t+b(1-\theta) \mathrm{d} W_t
\end{array}\right.
\end{equation*}
to get the following splitting method 
\begin{equation}\label{eq:splitting}
\begin{split}
	Q_{n+1} & =\left[b^{2} \theta(\theta-1)\tau+\left[a_{0}\tau\left[2(\theta-1) a_{-1}\tau\left[(\gamma-1) a_{2}\tau+y_n^{\frac{\gamma-1}{\theta-1}}\right]^{\frac{1}{1-\gamma}}\right]^{\frac{1}{2}}\right]^{2(1-\theta)}\right]^{\frac{1}{2}}, \\
	{ y_{n+1} } & =e^{-(\theta-1) a_{1} \tau} Q_{n+1}+b(1-\theta)  e^{-(\theta-1) a_{1}\tau} \Delta W_n,
\end{split}
\end{equation}
where $y_n \approx Y_{t_n}$ and $\Delta W_{n}:=W_{t_{n+1}}-W_{t_{n}}$,
$ y_0 = Y_0 $.

Although, as we see in  Section \ref{sec:numerics}, this scheme is easy to implement and has good properties, the direct analysis of \eqref{eq:splitting} is not straightforward and does not generalise (for example to cases where it is not feasible to solve the nonlinear terms exactly). 

Instead we choose to combine the splitting with a tamed Euler-Maruyama (see e.g. \cite{hutzenthaler2012}) for the nonlinear drift.
%
%
Combining the solutions of the OU process 
$dz=-\lambda z dt + b(1-\theta)dW_t$
and the tamed Euler-Maruyama approximation of 
$dv = F(v) dt$ we obtain the tamed-splitting method (\TSM) 
\begin{equation}\label{eq:tame-numerical-solution-without-jump}
	Y_{n+1} 
	= e^{-\lambda \tau}
	\left( Y_n + \tau \Ftau(Y_n)
	+ b (1 - \theta ) \Delta W_n \right),
\end{equation}
where $\Delta W_{n}:=W_{t_{n+1}}-W_{t_{n}}$ and the taming function is given by
$$  \Ftau(x)  := \frac{F(x)}{ 1 + \tau | F(x) |^2 }.$$
We note that key to the analysis below is the following inequality
\begin{equation}\label{eq:tame-contronal}
	| \Ftau(x) |^2
	=
	\frac{ | F(x) |^2}{ ( 1 + \tau | F(x) |^2 )^2}
	\leq
	\frac{ \tau^{-1} + | F(x) |^2}{ 1 + \tau | F(x) |^2 }
	= \tau^{-1}.	
\end{equation} 
Given $Y_n>0$ we see from \eqref{eq:tame-numerical-solution-without-jump} there is no guarantee that $Y_{n+1}$ is positive (see also Section \ref{sec:positivity}). We combine this scheme with a positivity-preserving backstop scheme as in \cite{KLM2020} for the Cox-Ingersoll-Ross process. We choose the backward Euler method
\begin{equation}\label{eq:BEM_method}
	\hat{Y}_{n+1 } = \hat{Y}_{ n}+(F( \hat{Y}_{ n + 1 })  -\lambda \hat{Y}_{n+1} ) \tau
	+ b (1 - \theta ) \Delta W_n
\end{equation}
where from \cite{lei2021first} we have 
\begin{equation*}
  \big\| Y_{t_{ n }} 
  - \hat{ Y }_{ n } \big\|_{ L^2 ( \Omega ; \R ) }
  \leq  C \tau.
\end{equation*}
We therefore define the tamed-splitting method with backstop:
\begin{equation}\label{eq:numerical-backstop}
  \begin{split}
    \bar{Y}_{n+1} 
    = & \left( e^{-\lambda \tau}
    ( \bar{Y}_n 
    + \tau \Ftau(\bar{Y}_n)
    + b (1 - \theta ) \Delta W_n )
	\right) \cdot
	\mathbbm{ 1 }_{ \{ \bar{Y}_{ n + 1 } > 0 \} } \\
	& +
	\left( \bar{Y}_n +( 
	F( \bar{Y}_{ n + 1 }  ) -\lambda \bar{Y}_{n+1}  )\tau
	+ b (1 - \theta ) \Delta W_n 
	\right) \cdot
	\mathbbm{ 1 }_{ \{ \bar{Y}_{ n + 1} < 0  \} }.
	\end{split}
\end{equation}
For our analysis we define a continuous version of the tamed-splitting method 
\eqref{eq:tame-numerical-solution-without-jump}
as
\begin{equation}\label{eq:continuous-numerical-solution}
	\bar{Y}_t = e^{-\lambda t} Y_{0} 
	+ \int_{0}^{t} 
	 e^{- \lambda ( t - \lfloor s \rfloor ) } 
	 \Ftau ( \bar{Y}_{\lfloor s \rfloor} ) 
	\dd s
	+ \int_{0}^{t} e^{- \lambda ( t - \lfloor s \rfloor ) }
	b (1 - \theta ) \dd W_s.
\end{equation}
We see in Section \ref{sec:numerics} that the splitting method \eqref{eq:splitting} and tamed-splitting method \eqref{eq:tame-numerical-solution-without-jump} have the same rate of convergence, however  \eqref{eq:splitting} is more accurate in terms of the error constant.


\section{Bounded moments and strong convergence of the tamed-splitting method}
\label{sec:convergence}
\subsection{Bounded moments of the tamed-splitting method}
Before proving the moment bounds 
of the numerical solution, we give a preliminary lemma.
\begin{lemma}\label{lem:preliminary_bound-moment-numerical}
For the approximation process $ \bar{Y}_t $ given by 
	\eqref{eq:continuous-numerical-solution}, we have that 
	\begin{equation*}
 	I : = \bar{Y}_{ \lfloor s \rfloor} \cdot
   \Ftau(\bar{Y}_{\lfloor s \rfloor}) 
   = \frac{\bar{Y}_{ \lfloor s \rfloor} \cdot
   F(\bar{Y}_{ \lfloor s \rfloor}) }{ 1 + \tau | F(\bar{Y}_{ \lfloor s \rfloor}) |^2 }< C \max{\{ \tau^{\frac{1}{2} }, 
   1, \bar{Y}_{ \lfloor s \rfloor} \} },
   \quad s \in [0,T].
  \end{equation*}
  \end{lemma}
  \textbf{Proof of Lemma \ref{lem:preliminary_bound-moment-numerical}} 
We divide into three cases:
\begin{itemize}
	\item [ $\cdot$ ]
	Case 1 : 
	$ 0 < \bar{Y}_{ \lfloor s \rfloor} \leq \tau $.\\
	According to \eqref{eq:tame-contronal},	we have
	\begin{equation}\label{eq:esti-I1}
		I \leq
		 \bar{Y}_{ \lfloor s \rfloor} \cdot
		\tau^{ - \frac{ 1 }{ 2 } }
		\leq \tau \cdot
		\tau^{  - \frac{ 1 }{ 2 } } 
		\leq 
		\tau^{  \frac{ 1 }{ 2 } } .
	\end{equation}
	\item [ $\cdot$ ]  
	Case 2 : $ \bar{Y}_{ \lfloor s \rfloor} \geq 1 $. \\
	Since $ a_2 , a_{ - 1 }, a_{0} > 0 $ 
     and $\gamma > \theta > 1$, it is easy to show that there is a constant $M$ such that 
     \begin{equation}\label{eq:esti-I2}
     		F (x) 
     		= (\theta-1)\big[ a_2 x^{\frac{\gamma-\theta}{1-\theta}}
     		- a_{-1} x^{\frac{\theta+1} {\theta-1}} 
     		+ a_0 x^{\frac{\theta}{\theta-1}} 
     		+ \frac{1}{2} b^2 \theta x^{-1} \big]
     		\leq M , \quad \text{for}
     		\quad  x \geq 1.
     \end{equation}
 Thus 
 \begin{equation*}
 	I \leq
 	\frac{ M \bar{Y}_{ \lfloor s \rfloor} 
 		 }
 	{ 1 + \tau | F(\bar{Y}_{\lfloor s \rfloor}) |^2}
 	\leq M \bar{Y}_{ \lfloor s \rfloor}.
 \end{equation*}
    \item [ $\cdot$ ]
    Case 3 :
    $ \tau < \bar{Y}_{ \lfloor s \rfloor} < 1 $. \\
    Because 
    $ \bar{Y}_{ \lfloor s \rfloor} $
     has upper and lower bounds, we have 
     $ \bar{Y}_{ \lfloor s \rfloor}^{ \eta }  \leq 1 $
     for any $ \eta \geq 0 $ and 
     $ \bar{Y}_{ \lfloor s \rfloor}^{ \eta }  \leq \tau^{\eta} $
     for any $ \eta < 0 $.
     Therefore, we directly calculate
    \begin{align*}
    		I  = &
    	\frac{\bar{Y}_{ \lfloor s \rfloor} \cdot
    		F(\bar{Y}_{\lfloor s \rfloor}) }
    	{ 1 + \tau | F(\bar{Y}_{\lfloor s \rfloor}) |^2} \\
    	= &
    	 \frac{
    	 	( \theta - 1 )
    	 	\big[ a_2 \bar{Y}_{ \lfloor s \rfloor}^{ 
    	 	\frac{\gamma + 1 - 2 \theta }
    	 	{ 1 - \theta } } - 
     	a_{ - 1 } \bar{Y}_{ \lfloor s \rfloor}^{ 
     	\frac{ 2 \theta }{ \theta - 1  } }
     + a_{ 0 } 
     \bar{Y}_{ \lfloor s \rfloor}^{ 
     \frac{ 2 \theta - 1 }
     { \theta - 1 } }
      + \frac{ 1 }{ 2 }
      b^2 \theta
    	 	\big] }
    	 { 1 + \tau ( \theta - 1 )^2 \Big| 
    	 	a_2 \bar{Y}_{ \lfloor s \rfloor}^{ 
    	 		\frac{\gamma  -  \theta }
    	 		{ 1 - \theta } } - 
    	 	a_{ - 1 } \bar{Y}_{ \lfloor s \rfloor}^{ 
    	 		\frac{  \theta + 1 }{ \theta - 1  } }
    	 	+ a_{ 0 } 
    	 	\bar{Y}_{ \lfloor s \rfloor}^{ 
    	 		\frac{  \theta  }
    	 		{ \theta - 1 } }
    	 	+ \frac{ 1 }{ 2 }
    	 	b^2 \theta \bar{Y}_{ \lfloor s \rfloor}^{-1}
    	 	\Big|^2 }		.
\end{align*}
For the numerator we have that
\begin{equation}\label{eq:numerator}
  - a_{ - 1 } \bar{Y}_{ \lfloor s \rfloor}^{ 
     			\frac{ 2 \theta }{ \theta - 1  } }
     		+ a_{ 0 } 
     		\bar{Y}_{ \lfloor s \rfloor}^{ 
     			\frac{ 2 \theta - 1 }
     			{ \theta - 1 } }
     		+ \frac{ 1 }{ 2 }
     		b^2 \theta 
     		\leq
     	a_{ 0 } + 	\frac{ 1 }{ 2 }
     		b^2 \theta =: K.
\end{equation}
For the denominator we expand the square and note
\begin{equation}\label{eq:denominator}
  - 2 a_{ -1 } a_0 
     	 	\bar{Y}_{ \lfloor s \rfloor}^{ 
     	 		\frac{2 \theta  + 1 }
     	 		{  \theta - 1 } }
      		- a_{-1} b^2 \theta
      		\bar{Y}_{ \lfloor s \rfloor}^{ 
      			\frac{ 2  }
      			{  \theta - 1 } }
      		+	a_0 b^2 \theta 
      		\bar{Y}_{ \lfloor s \rfloor}^{ 
      			\frac{ 1 }
      			{  \theta - 1 } }
      			\geq
      	- 2 a_{ -1 } a_0 	- a_{-1} b^2 \theta
      	=: \bar{a}.
\end{equation}
By \eqref{eq:denominator} with $ \tau < \bar{Y}_{ \lfloor s \rfloor} < 1 $ yields
that for the denominator we have 
\begin{align*}
\Big| & a_2 \bar{Y}_{ \lfloor s \rfloor}^{ 
    	 		\frac{\gamma  -  \theta }
    	 		{ 1 - \theta } } - 
    	 	a_{ - 1 } \bar{Y}_{ \lfloor s \rfloor}^{ 
    	 		\frac{  \theta + 1 }{ \theta - 1  } }
    	 	+ a_{ 0 } 
    	 	\bar{Y}_{ \lfloor s \rfloor}^{ 
    	 		\frac{  \theta  }
    	 		{ \theta - 1 } }
    	 	+ \frac{ 1 }{ 2 }
    	 	b^2 \theta \bar{Y}_{ \lfloor s \rfloor}^{-1}
    	 	\Big|^2  \\
&  \geq 
a_2^2 \bar{Y}_{ \lfloor s \rfloor}^{\frac{2\gamma  -  2\theta }{ 1 - \theta } } 
+ a_{ - 1 }^2 \tau^{\frac{  2\theta + 2 }{ \theta - 1  } }
+ a_{ 0 }^2 \tau^{\frac{  2\theta  }{ \theta - 1 } }
+ \frac{ 1 }{ 4 }b^4 \theta Y_{ \lfloor s \rfloor}^{-2}
+  2 a_2 a_0 \bar{Y}_{ \lfloor s \rfloor}^{\frac{\gamma  - 2 \theta }{ 1 - \theta } }
- 2 a_2 a_{-1} \bar{Y}_{ \lfloor s \rfloor}^{ \frac{\gamma - 1 - 2 \theta}{1 - \theta}}
+  a_2 b^2 \theta \bar{Y}_{ \lfloor s \rfloor}^{\frac{\gamma  - 1 }{ 1 - \theta } }
+ \bar{a}.
\end{align*}
As a consequence 
\begin{align*}
I \leq  K \big( 1  
       	+ \bar{Y}_{ \lfloor s \rfloor}^{ \frac{\gamma + 1 - 2 \theta }
       		{ 1 - \theta } }
       	\big)  
\Big[ 1 + \tau ( \theta - 1 )^2 &
       \Big(
       	a_2^2 \bar{Y}_{ \lfloor s \rfloor}^{\frac{2\gamma  -  2\theta }{ 1 - \theta } } 
       		+ a_{ - 1 }^2 \tau^{\frac{  2\theta + 2 }{ \theta - 1  } }
       	+ a_{ 0 }^2 \tau^{\frac{  2\theta  }{ \theta - 1 } }+ \frac{ 1 }{ 4 }b^4 \theta \bar{Y}_{ \lfloor s \rfloor}^{-2} \\
& +  2 a_2 a_0 \bar{Y}_{ \lfloor s \rfloor}^{ \frac{\gamma  - 2 \theta }{ 1 - \theta }}
  - 2 a_2 a_{-1} \bar{Y}_{ \lfloor s \rfloor}^{ \frac{\gamma - 1 - 2 \theta}{ 1 - \theta } }
+  a_2 b^2 \theta \bar{Y}_{ \lfloor s \rfloor}^{\frac{\gamma  - 1 }{1 - \theta } }
+ \bar{a} \Big)  \Big]^{-1}.
\end{align*}
We multiply the numerator and denominator above by 
     $ \bar{Y}_{ \lfloor s \rfloor}^{ 
       		\frac{\gamma + 1 - 2 \theta }
       		{  \theta - 1 } }  $      	
to get       	
\begin{align*}
I \leq 
       \big[K \big( 1
       	+ \bar{Y}_{ \lfloor s \rfloor}^{ \frac{\gamma + 1 - 2 \theta }
       		{  \theta - 1 } }
       	\big) \big]
       \Big[ & \bar{Y}_{ \lfloor s \rfloor}^{ 
       		\frac{\gamma + 1 - 2 \theta }
       		{  \theta - 1 } }  
       	 + \tau ( \theta - 1 )^2 \Big(
       	a_2^2 \bar{Y}_{ \lfloor s \rfloor}^{ 
       		\frac{\gamma  -  1 }
       		{ 1 - \theta } } 
       	+ \frac{ 1 }{ 4 }
       	b^4 \theta \bar{Y}_{ \lfloor s \rfloor}^{
       		\frac{\gamma + 3 -  4 \theta }
       	{  \theta - 1 } }
       	+  2 a_2 a_0 
       	\bar{Y}_{ \lfloor s \rfloor}^{ 
       		\frac{ 1 }
       		{  \theta - 1 } }
       	+  a_2 b^2 \theta \bar{Y}_{ \lfloor s \rfloor}^{-2} \\
       	&  + 
       		\bar{Y}_{ \lfloor s \rfloor}^{ \frac{\gamma + 1 - 2 \theta }
       			{  \theta - 1 } }
       		\Big(  a_{ - 1 }^2 \tau^{ 
       			\frac{  2\theta + 2 }{ \theta - 1  } }
       		+ a_{ 0 }^2 
       		\tau^{ 
       			\frac{  2\theta  }
       			{ \theta - 1 } }
       		- 2 a_2 a_{-1} \bar{Y}_{ \lfloor s \rfloor}^{ 
     			\frac{\gamma - 1 - 2 \theta }
     			{ 1 - \theta } }
     			- \bar{a} \Big)\Big) \Big]^{-1}.
\end{align*}
     If $ \bar{Y}_{ \lfloor s \rfloor} \geq 
      \varepsilon_0 > \tau $, 
      where $ \varepsilon_0 $ is a constant 
      independent of $\tau$,
      it's same as case 2.
      Otherwise, $ Y_{ \lfloor s \rfloor} \rightarrow \tau $,
      we have
           \begin{align*}
       	 I \stackrel{ \bar{Y}_{ \lfloor s \rfloor}
       	 	 \rightarrow \tau }{\longrightarrow} 
       	\big[ K \big( 1
       		+ \tau^{ \frac{\gamma + 1 - 2 \theta }
       			{  \theta - 1 } }
       		\big) \big]
       	\Big[ &  \tau^{ 
       			\frac{\gamma + 1 - 2 \theta }
       			{  \theta - 1 } }        		+  ( \theta - 1 )^2 \Big(
       		a_2^2 \tau^{ 
       			\frac{\gamma  -  \theta }
       			{ 1 - \theta } } 
       		+ \frac{ 1 }{ 4 }
       		b^4 \theta \tau^{
       			\frac{\gamma + 2 -  3 \theta }
       			{  \theta - 1 } }
       		+  2 a_2 a_0 
       		\tau^{ 
       			\frac{ \theta }
       			{  \theta - 1 } }
       		+  a_2 b^2 \theta 
       		\tau^{ 
       			-1 } \Big) \\ 
       			& +  ( \theta - 1 )^2
       			\Big(  a_{ - 1 }^2 \tau^{ 
       				\frac{  \gamma + \theta + 2 }{ \theta - 1  } }
       			+ a_{ 0 }^2 
       			\tau^{ 
       				\frac{ \gamma + \theta  }
       				{ \theta - 1 } }
       			- 2 a_2 a_{-1} \tau^{ 
     			\frac{\theta + 1 }
     			{  \theta - 1 } }
       			+ \bar{a} \tau^{ \frac{\gamma - \theta  }
       				{  \theta - 1 } }
       			 \Big) \big]^{-1}.
    \end{align*}
We see that for arbitrarily small 
$ \tau $, we have $I< \infty$.
\end{itemize}
We have completed the three cases and the Lemma is proved.
\qed
In the following error analysis, the moment bounds 
of the numerical solution are frequently used.
The next lemma indicates when $ q $-th positive moments of the solution to \eqref{eq:tame-numerical-solution-without-jump} are bounded.
The proof of Lemma \ref{lem:bound-moment-numerical} and
Lemma \ref{lem:inverse-bound-moment-numerical}
is inspired by \cite{kumar2019milstein}.
\begin{lemma}\label{lem:bound-moment-numerical}
	The approximation process $ Y_{t} $ produced by 
	\eqref{eq:continuous-numerical-solution}
	obeys, for any $ q \geq 4 $,
	\begin{equation}\label{eq:bound_moment_numerical_solution}
		\sup_{t \in [0 , T] } \E [ | \bar{Y}_t |^{q} ]
		\leq C ( 1 + | Y_0 |^q ).
	\end{equation}
\end{lemma}
\textbf{Proof of Lemma \ref{lem:bound-moment-numerical}}
Applying the
It\^o formula and 
\eqref{eq:continuous-numerical-solution}, 
for any $ l \in \R $, yields
\begin{equation}\label{eq:Ito-formula}
	\begin{split}
		( 1+ |\bar{Y}_t|^{2})^{ l }   
		=&
		( 1+|Y_0|^{2} )^{ l }
		+ 2 l \int_{0}^{t}
		( 1+|\bar{Y}_{s}|^{2})^{l -1} 
		\bar{Y}_{s} \cdot
		e^{ - \lambda ( t - \lfloor s \rfloor) }
		\Ftau(\bar{Y}_{\lfloor s \rfloor}) 
		\, \dd s \\
		& 
		+ l 
		\int_{0}^{t}
		(1+|\bar{Y}_{s}|^{2})^{ l -1}
		| e^{ - \lambda ( t - \lfloor s \rfloor) }
		b ( 1- \theta )  |^{2} 
		\, \dd s \\
		&
		+ l ( l - 1 )
		\int_{0}^{t}
		(1+|\bar{Y}_{s}|^{2})^{l - 2 }
		| \bar{Y}_{s} \cdot 
		e^{ - \lambda ( t - \lfloor s \rfloor) }
		b ( 1- \theta )
		|^{2}
		\, \dd s \\
		& 
		+ 2 l 
		\int_{0}^{t}
		( 1 + | \bar{Y}_{s}|^{2} )^{  l - 1 }
		Y_{s} \cdot
		e^{ - \lambda ( t - \lfloor s \rfloor) }
		b ( 1- \theta ) 
		\, \dd W_s .
	\end{split}
\end{equation}
Let $ l = \frac{ q }{ 2 } $, and taking expectation we get 
\begin{align*}
\E [ ( 1+  & |\bar{Y}_t|^{2})^{\frac{ q }{2}}  ] 
		\leq   
		( 1+|Y_0|^{2} )^{\frac{ q }{2}}
		+ q \E \left[  \int_{0}^{t}
		( 1+|\bar{Y}_{s}|^{2})^{\frac{ q }{2}-1} 
		\bar{Y}_{s} \cdot
		e^{ - \lambda ( t - \lfloor s \rfloor) }
		\Ftau(\bar{Y}_{\lfloor s \rfloor}) 
		\, \dd s \right] \\
		&
		+ \frac{ q ( q - 1 ) } { 2 } 
		\E \left[ \int_{0}^{t}
		(1+|\bar{Y}_{s}|^{2})^{\frac{ q } { 2 } - 1 } 
		| e^{ - \lambda ( t - \lfloor s \rfloor) }
		b ( 1- \theta )
		|^{2}
		\, \dd s \right] .
		\end{align*}
	Adding and subtracting terms we make a further decomposition
as follows	
\begin{align}\label{eq:esti-inverse-bound-moment}
\E [ & ( 1+  |\bar{Y}_t|^{2})^{\frac{ q }{2}}  ] \leq 
		( 1+|Y_0|^{2} )^{\frac{ q }{2}} \nonumber \\
	& 	+ q \E \left[  \int_{0}^{t}
		( 1+|\bar{Y}_{s}|^{2})^{\frac{ q }{2}-1} 
		( \bar{Y}_{s} - e^{- \lambda ( s - \lfloor s \rfloor)}
		\bar{Y}_{ \lfloor s \rfloor} ) \cdot
		e^{ - \lambda ( t - \lfloor s \rfloor) }
		\Ftau(\bar{Y}_{\lfloor s \rfloor}) 
		\, \dd s \right] \nonumber \\
		&
		+  q \E \bigg[  \int_{0}^{t}
		( 1+|\bar{Y}_{s}|^{2})^{\frac{ q }{2}-1} 
		\Big( 
		e^{- \lambda ( s - \lfloor s \rfloor)}
		\bar{Y}_{ \lfloor s \rfloor} \cdot
		e^{ - \lambda ( t - \lfloor s \rfloor) }
		\Ftau(\bar{Y}_{\lfloor s \rfloor})  
		+ \frac{   q - 1  } { 2 } 
		| e^{ - \lambda ( t - \lfloor s \rfloor) }
		b ( 1- \theta )
		|^{2}
		\Big)
		\, \dd s \bigg] \nonumber \\
		& =: 
		( 1+|Y_0|^{2} )^{\frac{ q }{2}}
		+ J_1 + J_2 .
\end{align}
In what follows, we bound the two terms $ J_1 $, $ J_2 $
separately. By the definition of 
$ Y_t $, we have
\begin{align*}
		J_1 = &
		q \E \left[  \int_{0}^{t}
		( 1 + |\bar{Y}_{s}|^{2})^{\frac{ q }{2}-1} 
		\int_{\lfloor s \rfloor}^{s} 
		e^{- \lambda ( s - \lfloor r \rfloor )}
		\Ftau ( \bar{Y}_{ \lfloor r \rfloor } ) \dd r
		\cdot
		e^{ - \lambda ( t - \lfloor s \rfloor) }
		\Ftau(\bar{Y}_{\lfloor s \rfloor}) 
		\, \dd s \right] \\
		& +
		q \E \left[  \int_{0}^{t}
		( 1 + |\bar{Y}_{s}|^{2})^{\frac{ q }{2}-1} 
		\int_{\lfloor s \rfloor}^{s} 
		e^{- \lambda ( s - \lfloor r \rfloor )}
		b ( 1 - \theta ) \dd W_r
		\cdot
		e^{ - \lambda ( t - \lfloor s \rfloor) }
		\Ftau(\bar{Y}_{\lfloor s \rfloor}) 
		\, \dd s \right].
		\end{align*}
Adding in and subtracting out $( 1 + |\bar{Y}_{\lfloor s \rfloor}|^{2})^{\frac{ q }{2}-1}$
\begin{align*}\label{eq:esti-J1}
	& J_1 	\leq  
		 q   \E \left[  \int_{0}^{t}
		( 1 + |\bar{Y}_{s}|^{2})^{\frac{ q }{2}-1} 
		\int_{\lfloor s \rfloor}^{s} 
		| \Ftau ( \bar{Y}_{ \lfloor r \rfloor } ) | \dd r
		\cdot
		| \Ftau(\bar{Y}_{\lfloor s \rfloor}) |
		\, \dd s \right] \\
		& +
		q \E \left[  \int_{0}^{t}
		( 1 + |\bar{Y}_{\lfloor s \rfloor}|^{2})^{\frac{ q }{2}-1} 
		\int_{\lfloor s \rfloor}^{s} 
		e^{- \lambda ( s - \lfloor r \rfloor )}
		b ( 1 - \theta ) \dd W_r
		\cdot
		e^{ - \lambda ( t - \lfloor s \rfloor) }
		\Ftau(\bar{Y}_{\lfloor s \rfloor}) 
		\, \dd s \right] \\
		& +
		q \E \left[  \int_{0}^{t}
		\big( ( 1 + |\bar{Y}_{ s }|^{2})^{\frac{ q }{2}-1}
		- ( 1 + |\bar{Y}_{\lfloor s \rfloor}|^{2})^{\frac{ q }{2}-1} \big)
		\int_{\lfloor s \rfloor}^{s} 
		e^{- \lambda ( s - \lfloor r \rfloor )}
		b ( 1 - \theta ) \dd W_r
		\cdot
		e^{ - \lambda ( t - \lfloor s \rfloor) }
		\Ftau(\bar{Y}_{\lfloor s \rfloor}) 
		\, \dd s \right].
\end{align*}
Using \eqref{eq:tame-contronal}
and the It\^o Isometry, we get
\begin{align*}
	& J_1 	\leq  
		 q   \E \left[  \int_{0}^{t}
		( 1 + |\bar{Y}_{s}|^{2})^{\frac{ q }
		{2}-1} 
		( s - \lfloor s \rfloor )
		\cdot \tau 
		\, \dd s \right] \\
		& +
		q \E \left[  \int_{0}^{t}
		\big( ( 1 + |\bar{Y}_{ s }|^{2})^{\frac{ q }{2}-1}
		- ( 1 + |\bar{Y}_{\lfloor s \rfloor}|^{2})^{\frac{ q }{2}-1} \big)
		\int_{\lfloor s \rfloor}^{s} 
		e^{- \lambda ( s - \lfloor r \rfloor )}
		b ( 1 - \theta ) \dd W_r
		\cdot
		e^{ - \lambda ( t - \lfloor s \rfloor) }
		\Ftau(\bar{Y}_{\lfloor s \rfloor}) 
		\, \dd s \right].
\end{align*}
With the aid of \eqref{eq:Ito-formula}
taking $ l = \frac{ q }{ 2 } - 1 $,
we can split 
$ J_1 $ into four additional terms:
\begin{align*}
J_1 \leq & C \E \left[  \int_{0}^{t}
		( 1 + |\bar{Y}_{s}|^{2})^{\frac{ q }{2}} 
		\right]  +
		q ( q - 2 )
		\E \Big[  \int_{0}^{t}
		\int_{\lfloor s \rfloor}^{s} 
		( 1 + |\bar{Y}_{ r }|^{2})^{\frac{ q }{2} - 2 }
		\bar{Y}_{ r } \cdot 
		e^{- \lambda ( s - \lfloor r \rfloor )}
		\Ftau ( \bar{Y}_{ \lfloor r \rfloor } ) \dd r \\
		& \times
		\int_{\lfloor s \rfloor}^{s} 
		e^{- \lambda ( s - \lfloor r \rfloor )}
		b ( 1 - \theta ) \dd W_r
		\cdot
		e^{ - \lambda ( t - \lfloor s \rfloor) }
		\Ftau(\bar{Y}_{\lfloor s \rfloor}) 
		\, \dd s \Big] \\
		& +
		q( q - 2 ) \E \Big[  \int_{0}^{t}
		\int_{\lfloor s \rfloor}^{s} 
		( 1 + |\bar{Y}_{ r }|^{2})^{\frac{ q }{2} - 2 }
		\bar{Y}_{r} \cdot
		b ( 1 - \theta ) 
		e^{- \lambda ( s - \lfloor r \rfloor )} \dd W_r \\
		& \times
		\int_{\lfloor s \rfloor}^{s} 
		e^{- \lambda ( s - \lfloor r \rfloor )}
		b ( 1 - \theta ) \dd W_r
		\cdot
		e^{ - \lambda ( t - \lfloor s \rfloor) }
		\Ftau(\bar{Y}_{\lfloor s \rfloor}) 
		\, \dd s \Big] \\
		& +
		\frac{q ( q - 2 ) ( q - 4 ) }{ 2 }
		\E \Big[  \int_{0}^{t}
		\int_{\lfloor s \rfloor}^{s} 
		( 1 + |\bar{Y}_{ r }|^{2})^{\frac{ q }{2} - 3 }
		| \bar{Y}_{r} 
		b ( 1 - \theta ) 
		e^{- \lambda ( s - \lfloor r \rfloor )} |^2 \dd r \\
		& \times 
		\int_{\lfloor s \rfloor}^{s} 
		e^{- \lambda ( s - \lfloor r \rfloor )}
		b ( 1 - \theta ) \dd W_r
		\cdot
		e^{ - \lambda ( t - \lfloor s \rfloor) }
		\Ftau(\bar{Y}_{\lfloor s \rfloor}) 
		\, \dd s \Big] \\
		& +
		\frac{ q ( q - 2 ) }{ 2 }
		\E \Big[  \int_{0}^{t}	
		\int_{\lfloor s \rfloor}^{s} 
		( 1 + |\bar{Y}_{ r }|^{2})^{\frac{ q }{2} - 2 }
		|  
		b ( 1 - \theta ) 
		e^{- \lambda ( s - \lfloor r \rfloor )} |^2 \dd r \\
		& \times 
		\int_{\lfloor s \rfloor}^{s} 
		e^{- \lambda ( s - \lfloor r \rfloor )}
		b ( 1 - \theta ) \dd W_r
		\cdot
		e^{ - \lambda ( t - \lfloor s \rfloor) }
		\Ftau(\bar{Y}_{\lfloor s \rfloor}) 
		\, \dd s \Big] \\
		=: &
		C \E \left[  \int_{0}^{t}
		( 1 + |\bar{Y}_{s}|^{2})^{\frac{ q }{2}} 
		\right] 
		+ J_{11} + J_{12} + J_{13} + J_{14}.
\end{align*}
Next we estimate 
$ J_{11} $-$ J_{14} $
term by term. Using 
\eqref{eq:tame-contronal}
leads to
\begin{equation*}
 \begin{split}
		J_{11}	
		\leq &
		q ( q - 2 )
		\E \Big[  \int_{0}^{t}
		\int_{\lfloor s \rfloor}^{s} 
		( 1 + |\bar{Y}_{ r }|^{2})^{\frac{ q }{2} - 2 }
		| \bar{Y}_{ r } | \cdot 
		| \Ftau ( \bar{Y}_{ \lfloor r \rfloor } ) | \dd r \\
		& \times
		\Big| \int_{\lfloor s \rfloor}^{s} 
		e^{- \lambda ( s - \lfloor r \rfloor )}
		b ( 1 - \theta ) \dd W_r \Big|
		\cdot
		| \Ftau(\bar{Y}_{\lfloor s \rfloor}) |
		\, \dd s \Big]  \\
		\leq &
		q ( q - 2 )
		\tau^{-1} \E \Big[  \int_{0}^{t}
		\int_{\lfloor s \rfloor}^{s} 
		( 1 + |\bar{Y}_{ r }|^{2})^{\frac{ q }{2} - 1 }
		\dd r 
		\Big| \int_{\lfloor s \rfloor}^{s} 
		e^{- \lambda ( s - \lfloor r \rfloor )}
		b ( 1 - \theta ) \dd W_r \Big| \dd s \Big] . 
		\end{split}
		\end{equation*}
	We now use Young's inequality and the H\"older inequality to obtain
	\begin{align*}
	J_{11}	\leq &
		C
		\tau^{-1} \E \Big[  \int_{0}^{t}
		\Big( \int_{\lfloor s \rfloor}^{s} 
		( 1 + |\bar{Y}_{ r }|^{2})^{\frac{ q }{2} - 1 }
		\dd r \Big)^{\frac{ q }{ q - 1}}
		+
		\Big| \int_{\lfloor s \rfloor}^{s} 
		e^{- \lambda ( s - \lfloor r \rfloor )}
		b ( 1 - \theta ) \dd W_r \Big|^{ q } \dd s \Big]  \\
		\leq &
		C
		\tau^{-1} \E \Big[  \int_{0}^{t}
		( s - \lfloor s \rfloor )^{ \frac{ 1 }{ q -1 } }
		\int_{\lfloor s \rfloor}^{s} 
		( 1 + |\bar{Y}_{ r }|^{2})^{ \frac{ q }{ 2 }
			 \cdot \frac{ q - 2}{ q - 1} }
		\dd r   \dd s \Big]  \\
		& +
	     C
		\tau^{-1} \E \Big[  \int_{0}^{t}
		\Big(
		\int_{\lfloor s \rfloor}^{s} 
		| e^{- \lambda ( s - \lfloor r \rfloor )}
		b ( 1 - \theta ) |^2 \dd r \Big)^{\frac{ q }{ 2 }}
		  \dd s \Big] . 
		  \end{align*}
	Now, using the  It\^o isometry yields
		  \begin{align}\label{eq:esti_split_J_11}
		J_{11} \leq &
		C \tau^{-1}
		\int_{0}^{t}
		\sup_{ 0 \leq r \leq s}
		\E [( 1 + |\bar{Y}_{ r }|^{2})^{ \frac{ q }{ 2 } }
		( s - \lfloor s \rfloor )^{\frac{ q }{ q - 1}}  ]\dd s  
		+
		 C \tau^{-1}
		\E \Big[  \int_{0}^{t}
		| b ( 1 - \theta ) |^2  
			( s - \lfloor s \rfloor )^{\frac{ q }{ 2 }} 
		 \dd s \Big]  \nonumber \\
		\leq &
		C \int_{0}^{t}
		\sup_{ 0 \leq r \leq s}
		\E [( 1 + |\bar{Y}_{ r }|^{2})^{ \frac{ q }{ 2 } } ]  \dd s .
	\end{align}
For $J_{12}$, using the generalized It\^o isometry
(see, e.g., 
\cite[Theorem 2.3.4]{zhang2017numerical} ) and
\eqref{eq:tame-contronal}
yields
\begin{align*}
		J_{12}
		= &
		q( q - 2 )  \int_{0}^{t}
		\int_{\lfloor s \rfloor}^{s} 
		\E \big[ 
		( 1 + |\bar{Y}_{ r }|^{2})^{\frac{ q }{2} - 2 }
		\bar{Y}_{r} \cdot
		| b ( 1 - \theta ) |^2
		e^{- 2 \lambda ( s - \lfloor r \rfloor )}  \nonumber \\
		& \times
		e^{- \lambda ( s - \lfloor r \rfloor )}
		b ( 1 - \theta ) 
		\cdot
		e^{ - \lambda ( t - \lfloor s \rfloor) }
		\Ftau(\bar{Y}_{\lfloor s \rfloor}) 
		\big] \dd r		
		\, \dd s  \nonumber \\
		\leq &
		C \tau^{- \frac{ 1 } { 2 } }
		\int_{0}^{t}
		\int_{\lfloor s \rfloor}^{s}
		\E \big[
		( 1 + |\bar{Y}_{ r }|^{2})^{\frac{ q }{2} - 1 } \big]
		\dd r \, \dd s \nonumber  \\
		\leq &
		C \tau^{- \frac{ 1 } { 2 } }
		\int_{0}^{t}
			\sup_{ 0 \leq r \leq s}
			\E \big[
		( 1 + |\bar{Y}_{ r }|^{2})^{\frac{ q }{2} } \big]
		( s - \lfloor s \rfloor )
		\, \dd s \nonumber  \\
		\leq &
		C \int_{0}^{t}
		\sup_{ 0 \leq r \leq s}
		\E [( 1 + |\bar{Y}_{ r }|^{2})^{ \frac{ q }{ 2 } } ]  \dd s .
	\end{align*}
For $ J_{13}$, using 
\eqref{eq:tame-contronal},
the Young inequality, 
we obtain that
\begin{align*}
		J_{13}
		\leq &
		C
		\E \Big[  \int_{0}^{t}
		\int_{\lfloor s \rfloor}^{s} 
		( 1 + |\bar{Y}_{ r }|^{2})^{\frac{ q }{2} - 3 }
		| \bar{Y}_{r} 
		b ( 1 - \theta )  |^2 \dd r 
		\Big| \int_{\lfloor s \rfloor}^{s} 
		e^{- \lambda ( s - \lfloor r \rfloor )}
		\dd W_r \Big|
		\cdot
		\big| \Ftau(\bar{Y}_{\lfloor s \rfloor}) \big|
		\, \dd s \Big] 	\\
		= &
		C \tau^{-\frac{1}{2}}
		\E \Big[  \int_{0}^{t}
		\int_{\lfloor s \rfloor}^{s} 
		( 1 + |\bar{Y}_{ r }|^{2})^{\frac{ q }{2} - 2 } \dd r  
		\Big| \int_{\lfloor s \rfloor}^{s} 
		e^{- \lambda ( s - \lfloor r \rfloor )}
		\dd W_r \Big|
		\, \dd s \Big] 	\\
		\leq &
		C \tau^{-\frac{1}{2}}
		\int_{0}^{t}
		\E \Big[ \Big( \int_{\lfloor s \rfloor}^{s} 
		( 1 + |\bar{Y}_{ r }|^{2})^{\frac{ q -4  }{2} } \dd r 
		\Big)^{ \frac{ q }{ q - 1} } \Big] 
		+
		\E \Big[ \Big| \int_{\lfloor s \rfloor}^{s} 
		e^{- \lambda ( s - \lfloor r \rfloor )}
		\dd W_r \Big|^{ q } \Big]
		\, \dd s.  	
		\end{align*}
	Now, using the It\^o isometry and
	the H\"older inequality obtains
\begin{align*}
		J_{13}\leq &
		C \tau^{-\frac{1}{2}}
		\int_{0}^{t}
		 \int_{\lfloor s \rfloor}^{s} 
		 \E \big[
		( 1 + |\bar{Y}_{ r }|^{2})^{\frac{ q  }{2} 
		\cdot \frac{ q - 4 }{ q - 1 } }
		\big]  \dd r 
		( s - \lfloor s \rfloor )^{ 
			\frac{ 1 }{ q - 1} }  \dd s +
		C \tau^{-\frac{1}{2}}
		\int_{0}^{t}
		\E \Big[ \Big( \int_{\lfloor s \rfloor}^{s} 
		| e^{- \lambda ( s - \lfloor r \rfloor )} |^2
		\dd r \Big)^{ \frac{ q }{ 2 } } \Big]
		\, \dd s  \nonumber	\\
		\leq &
		C \tau^{-\frac{1}{2}}
		\int_{0}^{t}
		\E \big[
		( 1 + |\bar{Y}_{ r }|^{2})^{\frac{ q  }{2} } 
		\big]  \dd r 
		( s - \lfloor s \rfloor )^{ 
			\frac{ q }{ q - 1} }  \dd s +
		C \tau^{-\frac{1}{2}}
		\int_{0}^{t}
     	 ( s - \lfloor s \rfloor )^{ \frac{ q }{ 2 } }
		\, \dd s  \nonumber	\\
		\leq &
		C \int_{0}^{t}
		\sup_{ 0 \leq r \leq s}
		\E [( 1 + |\bar{Y}_{ r }|^{2})^{ \frac{ q }{ 2 } } ]  \dd s .
\end{align*}
Similar to $J_{13}$, we obtain
\begin{equation}\label{eq:esti_split_J_14}
	J_{14} \leq
	C \int_{0}^{t}
	\sup_{ 0 \leq r \leq s}
	\E [( 1 + |\bar{Y}_{ r }|^{2})^{ \frac{ q }{ 2 } } ]  \dd s.
\end{equation}
Thus, putting \eqref{eq:esti_split_J_11}-\eqref{eq:esti_split_J_14} together we get  
\begin{equation*}
	J_{1} \leq
	C \int_{0}^{t}
	\sup_{ 0 \leq r \leq s}
	\E [( 1 + |\bar{Y}_{ r }|^{2})^{ \frac{ q }{ 2 } } ]  \dd s.
\end{equation*}
Returing to $J_2$ in \eqref{eq:esti-inverse-bound-moment}, with the help of Lemma \ref{lem:preliminary_bound-moment-numerical},
one can show that
\begin{align*}
		J_{2} \leq &
		 q  \E \left[  \int_{0}^{t}
		( 1+|\bar{Y}_{s}|^{2})^{\frac{ q }{2}-1}
		\frac{\bar{Y}_{ \lfloor s \rfloor} \cdot
			F(\bar{Y}_{\lfloor s \rfloor}) }
		{ 1 + \tau | F(\bar{Y}_{\lfloor s \rfloor}) |^2} 
		\, \dd s \right] 
		+ \frac{ |	b ( 1- \theta )
			|^{2} q ( q - 1 ) } { 2 } 
		\E \left[ \int_{0}^{t}
		(1+|\bar{Y}_{s}|^{2})^{\frac{ q } { 2 } - 1 } 
		\, \dd s \right] \\
		\leq &
		C \E \left[  \int_{0}^{t}
		( 1+|\bar{Y}_{s}|^{2})^{\frac{ q }{2}-1}
		\big( 1 + | \bar{Y}_{ \lfloor s \rfloor} |^2
		\big)
		\, \dd s \right] 
		+ C
		\E \left[ \int_{0}^{t}
		(1+|\bar{Y}_{s}|^{2})^{\frac{ q } { 2 } - 1 } 
		\, \dd s \right] \\
		\leq &
		C \int_{0}^{t}
		\sup_{ 0 \leq r \leq s}
		\E [( 1 + |\bar{Y}_{ r }|^{2})^{ \frac{ q }{ 2 } } ] \dd s.
\end{align*}
Taking  the estimates of $ J_{1}$ and $ J_{2}$ into account, 
we derive from 
\eqref{eq:esti-inverse-bound-moment} 
that for any $ t \in [ 0, T ] $
\begin{equation*}
	\sup_{ 0 \leq s \leq t}
	\E [( 1 + |\bar{Y}_{ s }|^{2})^{ \frac{ q }{ 2 } } ]
	\leq 
	( 1 + |Y_{ 0 }|^{2})^{ \frac{ q }{ 2 } }
	+ C \int_{0}^{t}
	\sup_{ 0 \leq r \leq s}
	\E [( 1 + |\bar{Y}_{ r }|^{2})^{ \frac{ q }{ 2 } } ] \dd s.
\end{equation*}
 Applying the Gronwall inequality thus completes the proof.
\qed

Next, we need to prove inverse bound moments of
numerical solution.
\begin{lemma}\label{lem:inverse-bound-moment-numerical}
		The approximation process $ Y_t $ produced by 
	\eqref{eq:continuous-numerical-solution}
	obeys, for any $ q \leq -1 $,
	\begin{equation*}\label{eq:inverse_bound_moment_numerical_solution}
		\sup_{t \in [0 , T] } \E [ | \bar{Y}_t |^{q} ]
		\leq C ( 1 + | Y_0 |^q ).
	\end{equation*}	
\end{lemma}

\textbf{Proof of Lemma \ref{lem:inverse-bound-moment-numerical} }
By \eqref{eq:Ito-formula} with
$ l = \frac{ q }{ 2 } $,
$ q \leq -1 $,
similar to the proof of Lemma
\ref{lem:bound-moment-numerical},
we have
 \begin{align}\label{eq:esti_inverse_moment}
 	 \E [ ( 1+ & |\bar{Y}_t|^{2})^{\frac{ q }{2}}  ]
 		\leq 
 		( 1+|Y_0|^{2} )^{\frac{ q }{2}}
 		+
 		C \E \left[  \int_{0}^{t}
 		( 1 + |\bar{Y}_{s}|^{2})^{\frac{ q }{2}} 
 		\right] \nonumber \\
 		& +
 		q ( q - 2 )
 		 \E \Big[  \int_{0}^{t}
 		 	\int_{\lfloor s \rfloor}^{s} 
 		( 1 + |\bar{Y}_{ r }|^{2})^{\frac{ q }{2} - 2 }
 		Y_{r } \cdot 
 		e^{- \lambda ( s - \lfloor r \rfloor )}
 		\Ftau ( \bar{Y}_{ \lfloor r \rfloor } ) \dd r \nonumber \\
 		 & \times
 		\int_{\lfloor s \rfloor}^{s} 
 		e^{- \lambda ( s - \lfloor r \rfloor )}
 		b ( 1 - \theta ) \dd W_r
 		\cdot
 		e^{ - \lambda ( t - \lfloor s \rfloor) }
 		\Ftau(\bar{Y}_{\lfloor s \rfloor}) 
 		\, \dd s \Big] \nonumber \\
 		& +
 		q( q - 2 ) \E \Big[  \int_{0}^{t}
 		 \int_{\lfloor s \rfloor}^{s} 
 		( 1 + |\bar{Y}_{ r }|^{2})^{\frac{ q }{2} - 2 }
 		\bar{Y}_{r} \cdot
 		b ( 1 - \theta ) 
 		e^{- \lambda ( s - \lfloor r \rfloor )} \dd W_r \nonumber \\
 		 & \times
 		\int_{\lfloor s \rfloor}^{s} 
 		e^{- \lambda ( s - \lfloor r \rfloor )}
 		b ( 1 - \theta ) \dd W_r
 		\cdot
 		e^{ - \lambda ( t - \lfloor s \rfloor) }
 		\Ftau(\bar{Y}_{\lfloor s \rfloor}) 
 		\, \dd s \Big] \nonumber \\
 		& +
 		 \frac{q ( q - 2 ) ( q - 4 ) }{ 2 }
 		\E \Big[  \int_{0}^{t}
 		\int_{\lfloor s \rfloor}^{s} 
 		( 1 + |\bar{Y}_{ r }|^{2})^{\frac{ q }{2} - 3 }
 		| \bar{Y}_{r} 
 		b ( 1 - \theta ) 
 		e^{- \lambda ( s - \lfloor r \rfloor )} |^2 \dd r \nonumber  \\
 		 & \times 
 		\int_{\lfloor s \rfloor}^{s} 
 		e^{- \lambda ( s - \lfloor r \rfloor )}
 		b ( 1 - \theta ) \dd W_r
 		\cdot
 		e^{ - \lambda ( t - \lfloor s \rfloor) }
 		\Ftau(\bar{Y}_{\lfloor s \rfloor}) 
 		\, \dd s \Big] \nonumber  \\
 		& +
 		 \frac{ q ( q - 2 ) }{ 2 }
 		 \E \Big[  \int_{0}^{t}	
 		\int_{\lfloor s \rfloor}^{s} 
 		( 1 + |\bar{Y}_{ r }|^{2})^{\frac{ q }{2} - 2 }
 		|  
 		b ( 1 - \theta ) 
 		e^{- \lambda ( s - \lfloor r \rfloor )} |^2 \dd r \nonumber \\
 		& \times 
 		\int_{\lfloor s \rfloor}^{s} 
 		e^{- \lambda ( s - \lfloor r \rfloor )}
 		b ( 1 - \theta ) \dd W_r
 		\cdot
 		e^{ - \lambda ( t - \lfloor s \rfloor) }
 		\Ftau(\bar{Y}_{\lfloor s \rfloor}) 
 		\, \dd s \Big] \nonumber  \\
 		& +  
 		q \E \Big[  \int_{0}^{t}
 		( 1+|\bar{Y}_{s}|^{2})^{\frac{ q }{2}-1} 
 		\Big( 
 		e^{- \lambda ( s - \lfloor s \rfloor)}
 		\bar{Y}_{ \lfloor s \rfloor} \cdot
 		e^{ - \lambda ( t - \lfloor s \rfloor) }
 		\Ftau(\bar{Y}_{\lfloor s \rfloor}) 
 		\nonumber  \\
 		&
 		+ \frac{ ( q - 1 ) } { 2 } 
 		| e^{ - \lambda ( t - \lfloor s \rfloor) }
 		b ( 1- \theta )
 		|^{2} \Big)
 		\, \dd s \Big] \nonumber  \\
 		=: &
 		( 1+|Y_0|^{2} )^{\frac{ q }{2}}
 		 +
 		C \E \left[  \int_{0}^{t}
 		( 1 + |\bar{Y}_{s}|^{2})^{\frac{ q }{2}} 
 		\right] 
 		+ I_{1} + I_{2} + I_{3} + I_{4} +I_{5}.
 \end{align}
Next we estimate 
$ I_{1} $-$ I_{5} $
one by one. Using 
\eqref{eq:tame-contronal}
and the Young
inequality leads to
\begin{align*}
	I_{1}	
	\leq &
		q ( q - 2 )
	\E \Big[  \int_{0}^{t}
	\int_{\lfloor s \rfloor}^{s} 
	( 1 + |\bar{Y}_{ r }|^{2})^{\frac{ q }{2} - 2 }
	| \bar{Y}_{ r } | \cdot 
	| \Ftau ( \bar{Y}_{ \lfloor r \rfloor } ) | \dd r \\
	& \times
	\Big| \int_{\lfloor s \rfloor}^{s} 
	e^{- \lambda ( s - \lfloor r \rfloor )}
	b ( 1 - \theta ) \dd W_r \Big|
	\cdot
	| \Ftau(\bar{Y}_{\lfloor s \rfloor}) |
	\, \dd s \Big]  \\
	\leq &
	q ( q - 2 )
	\tau^{-1} \E \Big[  \int_{0}^{t}
	\int_{\lfloor s \rfloor}^{s} 
	( 1 + |\bar{Y}_{ r }|^{2})^{\frac{ q }{2} - 1 }
	 \dd r 
	\Big| \int_{\lfloor s \rfloor}^{s} 
	e^{- \lambda ( s - \lfloor r \rfloor )}
	b ( 1 - \theta ) \dd W_r \Big| \dd s \Big]  \\
		\leq &
	\frac{q ( q - 2 )} {2}
	\tau^{-1} \E \Big[  \int_{0}^{t}
	\Big( \int_{\lfloor s \rfloor}^{s} 
	( 1 + |\bar{Y}_{ r }|^{2})^{\frac{ q }{2} - 1 }
	\dd r \Big)^2 
	+
	\Big| \int_{\lfloor s \rfloor}^{s} 
	e^{- \lambda ( s - \lfloor r \rfloor )}
	b ( 1 - \theta ) \dd W_r \Big|^2 \dd s \Big]  
\end{align*}
With the help of the H\"older inequality and the 
It\^o isometry, one can get
\begin{align}\label{eq:esti_I1}
    I_{1}	\leq &
	\frac{q ( q - 2 )} {2}
	\tau^{-1} \E \Big[  \int_{0}^{t}
	( s - \lfloor s \rfloor )
	 \int_{\lfloor s \rfloor}^{s} 
	( 1 + |\bar{Y}_{ r }|^{2})^{ q  - 2 }
	\dd r   \dd s \Big] \nonumber  \\
	 & +
	\frac{q ( q - 2 )} {2}
	\tau^{-1} \E \Big[  \int_{0}^{t}
	\int_{\lfloor s \rfloor}^{s} 
	| e^{- \lambda ( s - \lfloor r \rfloor )}
	b ( 1 - \theta ) |^2 \dd r  \dd s \Big]  \nonumber \\
		\leq &
	\frac{q ( q - 2 )} {2}
	 \int_{0}^{t}
	\sup_{ 0 \leq r \leq s}
	\E [( 1 + |\bar{Y}_{ r }|^{2})^{ \frac{ q }{ 2 } }
	 ( s - \lfloor s \rfloor )  ]\dd s  
	 +
	\frac{q ( q - 2 )} {2}
	\E \Big[  \int_{0}^{t}
	| b ( 1 - \theta ) |^2   \dd s \Big]  \nonumber \\
	\leq &
	C \int_{0}^{t}
	\sup_{ 0 \leq r \leq s}
	\E [( 1 + |\bar{Y}_{ r }|^{2})^{ \frac{ q }{ 2 } } ]  \dd s .
  \end{align}
Also, using the Young inequality and the It\^o isometry yields
\begin{equation*}
	\begin{split}
		I_{2}
		\leq &
		\frac{q( q - 2 )} { 2 } 
		 \int_{0}^{t}
		\E \Big[ \Big| \int_{\lfloor s \rfloor}^{s} 
		( 1 + |\bar{Y}_{ r }|^{2})^{\frac{ q }{2} - 2 }
		Y_{r} \cdot
		b ( 1 - \theta ) 
		e^{- \lambda ( s - \lfloor r \rfloor )} 
		\dd W_r \Big|^2 \Big]
		\, \dd s  \\
		& +
		\frac{q( q - 2 )} { 2 } 
		  \int_{0}^{t}
		\E \Big[ \Big| \int_{\lfloor s \rfloor}^{s} 
		e^{- \lambda ( s - \lfloor r \rfloor )}
		b ( 1 - \theta ) 
		\cdot
		e^{ - \lambda ( t - \lfloor s \rfloor) }
		\Ftau(\bar{Y}_{\lfloor s \rfloor}) 
		\dd W_r \Big|^2 \Big]
		\, \dd s  \\
		\leq &
		\frac{q( q - 2 )} { 2 } 
		  \int_{0}^{t}
		 \E \Big[ \int_{\lfloor s \rfloor}^{s} 
		( 1 + |\bar{Y}_{ r }|^{2})^{q - 4 }
		| \bar{Y}_{r} |^2 \cdot
		| b ( 1 - \theta ) |^2
		\dd r  \Big]
		\, \dd s  \\
		& +
		\frac{q( q - 2 )} { 2 } \tau^{-1}
		 \int_{0}^{t}
		\E \Big[  \int_{\lfloor s \rfloor}^{s} 
		| b ( 1 - \theta ) |^2
		\dd r \Big]
		\, \dd s  .
	\end{split}
\end{equation*}
Because $ q \leq -1 $, it's easy to see that
$( 1 + |\bar{Y}_{ r }|^{2})^{ \frac{ q }{ 2 }}
\geq ( 1 + |\bar{Y}_{ r }|^{2})^{q - 3 }
\geq ( 1 + |\bar{Y}_{ r }|^{2})^{q - 4 }
\cdot |\bar{Y}_{ r }|^{2} .$
Thus
\begin{equation*}
    \begin{split}
    I_{2} \leq &
		\frac{q( q - 2 )} { 2 } 
		  \int_{0}^{t}
		  \E \Big[
		\int_{\lfloor s \rfloor}^{s} 
		( 1 + |\bar{Y}_{ r }|^{2})^{ \frac{ q }{ 2 }}
		| b ( 1 - \theta ) |^2
		\dd r  \Big]
		\, \dd s  
		 +
		\frac{q( q - 2 )} { 2 } 
		  \int_{0}^{t} 
		| b ( 1 - \theta ) |^2
		\, \dd s  \\
		\leq &
		C \int_{0}^{t}
		\sup_{ 0 \leq r \leq s}
		\E [( 1 + |\bar{Y}_{ r }|^{2})^{ \frac{ q }{ 2 } } ]  \dd s .
    \end{split}
\end{equation*}
For $I_3$, using \eqref{eq:tame-contronal}, Young's inequality
and the It\^o isometry,
we obtain that
\begin{align*}
	I_{3}
	\leq &
	\frac{q b ( q - 2 ) ( q - 4 ) (1 - \theta) }{ 2 }
	\E \Big[  \int_{0}^{t}
	\int_{\lfloor s \rfloor}^{s} 
	( 1 + |\bar{Y}_{ r }|^{2})^{\frac{ q }{2} - 3 }
	| \bar{Y}_{r} 
	b ( 1 - \theta )  |^2 \dd r \nonumber \\
	& \times 
	\Big| \int_{\lfloor s \rfloor}^{s} 
	e^{- \lambda ( s - \lfloor r \rfloor )}
	 \dd W_r \Big|
	\cdot
	\big| \Ftau(\bar{Y}_{\lfloor s \rfloor}) \big|
	\, \dd s \Big] 	\nonumber \\
		\leq &
	C \tau^{-\frac{1}{2}}
	\E \Big[  \int_{0}^{t}
	\int_{\lfloor s \rfloor}^{s} 
	( 1 + |\bar{Y}_{ r }|^{2})^{\frac{ q }{2} - 2 } \dd r  
	\Big| \int_{\lfloor s \rfloor}^{s} 
	e^{- \lambda ( s - \lfloor r \rfloor )}
	\dd W_r \Big|
	\, \dd s \Big] 	\nonumber \\
	\leq &
	C \tau^{-\frac{1}{2}}
	  \int_{0}^{t}
	\E \Big[ \Big( \int_{\lfloor s \rfloor}^{s} 
	( 1 + |\bar{Y}_{ r }|^{2})^{\frac{ q }{2} - 2 } \dd r 
	\Big)^2 \Big] \dd s
	+
		C \tau^{-\frac{1}{2}}
	\int_{0}^{t} 
	\E \Big[ \Big| \int_{\lfloor s \rfloor}^{s} 
	e^{- \lambda ( s - \lfloor r \rfloor )}
	\dd W_r \Big|^2 \Big]
	\, \dd s  	\nonumber \\
	\leq &
	C \int_{0}^{t}
	\sup_{ 0 \leq r \leq s}
	\E [( 1 + |\bar{Y}_{ r }|^{2})^{ \frac{ q }{ 2 } } ]  \dd s .
	\end{align*}
Similar to $ I_{3} $, we can obtain
\begin{equation*}
	I_{4} \leq
	C \int_{0}^{t}
	\sup_{ 0 \leq r \leq s}
	\E [( 1 + |\bar{Y}_{ r }|^{2})^{ \frac{ q }{ 2 } } ]  \dd s.
\end{equation*}
Since $ q \leq -1 $, we have 
\begin{align*}
    	I_{5} = &
		 q e^{- \lambda ( t + s -2\lfloor s \rfloor )}
		  \E \left[  \int_{0}^{t}
		( 1+|\bar{Y}_{s}|^{2})^{\frac{ q }{2}-1} 
		\bar{Y}_{ \lfloor s \rfloor} \cdot
		\Ftau(\bar{Y}_{\lfloor s \rfloor}) 
		\, \dd s \right] \\
		&
		+ e^{- 2\lambda ( t  - \lfloor s \rfloor )}
		\frac{ q ( q - 1 ) } { 2 } 
		\E \left[ \int_{0}^{t}
		(1+|\bar{Y}_{s}|^{2})^{\frac{ q } { 2 } - 1 } 
		|
		b ( 1- \theta )
		|^{2}
		\, \dd s \right] \\
		= &
		 q e^{- \lambda ( t+s-\lfloor s \rfloor )}
		 \E \left[  \int_{0}^{t}
		( 1+|\bar{Y}_{s}|^{2})^{\frac{ q }{2}-1}
		\frac{\bar{Y}_{ \lfloor s \rfloor} \cdot
		 F(\bar{Y}_{\lfloor s \rfloor}) }
		{ 1 + \tau | F(\bar{Y}_{\lfloor s \rfloor}) |} 
		\, \dd s \right] \\
		&
		+ \frac{e^{- 2\lambda ( t  - \lfloor s \rfloor )} |	b ( 1- \theta )
			|^{2} q ( q - 1 ) } { 2 } 
		\E \left[ \int_{0}^{t}
		(1+|\bar{Y}_{s}|^{2})^{\frac{ q } { 2 } - 1 } 
		\, \dd s \right] \\
		= &
     	 | q | e^{- \lambda ( t+s-\lfloor s \rfloor )} 
		 \int_{0}^{t} \E \Big[
		 	\frac{ 1  }
		 { 1 + \tau | F(\bar{Y}_{\lfloor s \rfloor}) | }
		 ( 1+|\bar{Y}_{s}|^{2})^{\frac{ q }{2}-1} 
		\\
		& \times 
		( \theta - 1 ) 
		\big[ - a_2 \bar{Y}_{\lfloor s \rfloor}^{
			\frac{ \gamma + 1 - 2 \theta}
			 {1 - \theta } } + a_{-1}
		 \bar{Y}_{\lfloor s \rfloor}^{
		 	\frac{ 2 \theta} 
		 	{ \theta - 1 } } 
	 	- a_0 \bar{Y}_{\lfloor s \rfloor}^{
	 		\frac{  2 \theta - 1 }
	 		{ \theta -1  } }
 		- \frac{1}{2} b^2  \theta ( \theta - 1 ) \big]
		\, \dd s \Big] \\
		&
		+ \frac{e^{- 2\lambda ( t  - \lfloor s \rfloor )} |	b ( 1- \theta )
			|^{2} q ( q - 1 ) } { 2 } 
		 \int_{0}^{t} \E \left[
		(1+|\bar{Y}_{s}|^{2})^{\frac{ q } { 2 } - 1 } 
		 \right]
		\, \dd s .
\end{align*}
This, together with the Young inequality and Lemma \ref{lem:bound-moment-numerical}, leads to
\begin{align}\label{eq:esti_I5}
		I_{5}\leq &
		| q | ( \theta - 1  )
		\int_{0}^{t} \E \big[
		( 1+|Y_{s}|^{2})^{	\frac{ q } { 2 } - 1 } 
		\cdot a_{-1}
		Y_{\lfloor s \rfloor}^{
			\frac{ 2 \theta} 
			{ \theta - 1 } } 
		\big] \, \dd s
		+ C
		\int_{0}^{t} \E \left[
		(1+|Y_{s}|^{2})^{\frac{ q } { 2 } - 1 } 
		\right]
		\, \dd s \nonumber \\
		\leq &
		C
		\int_{0}^{t} \E \big[
		( 1+|Y_{s}|^{2})^{	 q   - 2 } \big]
		+ \E \big[
		| Y_{\lfloor s \rfloor}|^{
			\frac{ 4 \theta} { \theta - 1 } } 
		\big] \, \dd s
		+ C
		\int_{0}^{t} \E \left[
		(1+|Y_{s}|^{2})^{\frac{ q } { 2 } - 1 } 
		\right]
		\, \dd s \nonumber \\
		\leq &
		C \int_{0}^{t}
		\sup_{ 0 \leq r \leq s}
		\E [( 1 + |Y_{ r }|^{2})^{ \frac{ q }{ 2 } } ] \dd s.
	\end{align}
Taking \eqref{eq:esti_I1}-\eqref{eq:esti_I5}
into account, 
we derive from 
\eqref{eq:esti_inverse_moment} 
that for any  $ t \in [ 0, T ] $. 
\begin{equation*}
	\sup_{ 0 \leq s \leq t}
	\E [( 1 + |Y_{ s }|^{2})^{ \frac{ q }{ 2 } } ]
	\leq 
	( 1 + |Y_{ 0 }|^{2})^{ \frac{ q }{ 2 } }
	+ C \int_{0}^{t}
	\sup_{ 0 \leq r \leq s}
	\E [( 1 + |Y_{ r }|^{2})^{ \frac{ q }{ 2 } } ] \dd s.
\end{equation*}
Applying the Gronwall inequality thus completes the proof.
\qed
\begin{remark}
Lemma \ref{lem:bound-moment-numerical} proves the bounded moment of the numerical solution
for $ q \geq 4 $. 
Infact, by applying the H\"older inequality, 
we can extend the results to $ 1 \leq q < 4 $.	
\end{remark}

\subsection{Strong convergence}
Before giving the convergence result, we introduce the following auxiliary process,
\begin{equation}\label{eq:numerical-solution-without-jump}
	\widetilde{Y}_{n+1} = e^{-\lambda \tau}
	( \widetilde{Y}_n +  
	\tau F ( Y_{t_{ n }} )
	+ b (1 - \theta ) \Delta W_n ),
\end{equation}
where $\lambda:=(\theta-1)a_1$ (as in \eqref{eq:defLamF}).
Then, we give the bounded moment of auxiliary process
$ \widetilde{Y}_{n} $.
\begin{lemma}\label{lem:bounded_moment_auxiliary}
For the auxiliary process defined in \eqref{eq:numerical-solution-without-jump}, we have
\begin{equation*}\label{eq:bounded_moment_auxiliary}
	\E [ |\widetilde{Y}_{n} |^2 ] < \infty,
	\quad
	\text{for} \quad n = 1, 2, \cdots, N.
\end{equation*}
\end{lemma}
\textbf{ Proof of Lemma \ref{lem:bounded_moment_auxiliary}}
It is easy to check that
\begin{equation*}
	\big\| \widetilde{Y}_{1} \big\|^{ 2  
	}_{ L^{2} ( \Omega ; \R )} 
	=
	| e^{-\lambda \tau } |^2
	\big\| Y_{0} + F ( Y_{t_0} ) \tau
	 + b ( 1 - \theta ) \Delta W_0 \big\|^{ 2  
	 }_{ L^{2} ( \Omega ; \R )}  
	 < \infty.	
\end{equation*}
Assume 
$ 	\big\| \widetilde{Y}_{ n - 1 } \big\|^{ 2  
}_{ L^{2 } ( \Omega ; \R )}  < \infty  $,
then by the definition of 
$ \widetilde{Y}_{n} $,
we have
\begin{align*}
		| \widetilde{Y}_{n} \big|^{ 2 }  
  =	&
  | e^{-\lambda \tau } |^2
  \big| \widetilde{Y}_{n-1} +
       F ( Y_{t_{n-1}}) \tau  
       + b ( \theta -1) \Delta W_{n-1} \big|^{ 2 }  \\
  \leq &
   \big| \widetilde{Y}_{n-1} \big|^{ 2 } 
    + \big| F ( Y_{t_{n-1}} ) \big|^{ 2 }  \tau^2
    + b^2 ( \theta - 1)^2 | \Delta W_{n-1} |^2
    + 2 \left\langle \widetilde{Y}_{n-1},
        F ( Y_{t_{n-1}} ) \right\rangle  \\
    & 
    + 2 \left\langle \widetilde{Y}_{n-1},
   b ( \theta -1) \Delta W_{n-1}  \right\rangle 
   + 2 \left\langle F ( Y_{t_{n-1}} ) \tau ,
   b ( \theta -1) \Delta W_{n-1}  \right\rangle .	
\end{align*}
Take the expectation of both sides, by the bounded and inverse bounded moment of $Y_t$, we have
\begin{align*}
		\big\| \widetilde{Y}_{n} \big\|^{ 2  
		}_{ L^{2 } ( \Omega ; \R )}  
		\leq 
		&
		\big\| \widetilde{Y}_{n-1} \big\|^{ 2  
		}_{ L^{2  } ( \Omega ; \R )} 
		+ \big\| F ( Y_{t_{n-1}} ) \big\|^{ 2  
		}_{ L^{2  } ( \Omega ; \R )}  \tau^2
		+ b^2 ( \theta - 1)^2 \| \Delta W_{n-1} \big\|^{ 2  
		}_{ L^{2 } ( \Omega ; \R )}  \\
	    &
		+ 2 \E \Big[ \left\langle \widetilde{Y}_{n-1},
		F ( Y_{t_{n-1}} ) \right\rangle 
		\Big] 
		+ 2 \E \Big[ \left\langle \widetilde{Y}_{n-1},
		b ( \theta -1) \Delta W_{n-1}  \right\rangle 
		\Big]
		 \\
		&
		+ 2 \E \Big[ \left\langle F ( Y_{t_{n-1}} ) \tau ,
		b ( \theta -1) \Delta W_{n-1}  \right\rangle
		\Big] \\
			\leq 
		&
		3 \big\| \widetilde{Y}_{n-1} \big\|^{ 2  
		}_{ L^{2  } ( \Omega ; \R )} 
		+ 
		3 \big\| F ( Y_{t_{n-1}} ) \big\|^{ 2  
		}_{ L^{2  } ( \Omega ; \R )}  \tau^2
		+ 
		3b^2 ( \theta - 1)^2 \tau^{  2  } 
		< \infty.
\end{align*}
This completes the lemma. \qed
We are now ready to prove the mean-square convergence rate of order one for the scheme.

\begin{theorem}\label{thm:convergence-rate-without-jump}
	Let $\{ Y_t \}_{0 \leq t \leq T }$
	and $\{ Y_{n} \}_{0 \leq n \leq N } $ be solutions to \eqref{eq:lamperti-translation1-jump} 
	and \eqref{eq:tame-numerical-solution-without-jump}, respectively. 
	Then for any time step $\tau>0$ it holds that
	\begin{equation*}\label{eq:convergence-rate}
		\sup _{N \in \mathbb{N}} \sup _{0 \leq n \leq N}
		\|Y_{n} -  Y_{t_{n}} \|_{L^{2}
			(\Omega ; \mathbb{R})} \leq C \tau .
	\end{equation*}
\end{theorem}
\textbf{Proof of Theorem \ref{thm:convergence-rate-without-jump}}
For $ t \in  [ t_{ n }, t_{ n + 1 } ] $,
we re-write the exact and numerical solution to
\eqref{eq:Lamperti},
\begin{equation}\label{eq:rewrite-exact-solution}
	Y_{t}
	=
	e^{ - \lambda ( t - t_{ n } ) } Y_{ t_{ n } }
	+
	\int_{ t_{n} }^{ t }
	e^{ - \lambda ( t - r ) }
	F ( Y_{ r } ) \dd r
	+
	\int_{t_{ n } }^{ t }
	e^{ - \lambda ( t - r ) }
	b ( 1 - \theta ) \dd W_{ r }
\end{equation}
and
\begin{equation}\label{eq:rewrite-numerical-solution}
	\bar{Y}_{t}
	=
	e^{ - \lambda ( t - t_{ n } ) } 
	\bar{Y}_{ t_{ n } }
	+
	\int_{ t_{n} }^{ t }
	e^{ - \lambda ( t - t_{n} ) }
	F^{ ( \tau ) } ( \bar{Y}_{ t_{ n } } ) \dd r
	+
	\int_{t_{ n } }^{ t }
	e^{ - \lambda ( t - t_{n} ) }
	b ( 1 - \theta ) \dd W_{ r }.
\end{equation}
We subtract \eqref{eq:rewrite-exact-solution}
from equation \eqref{eq:rewrite-numerical-solution}
to get
\begin{equation*}\label{eq:error}
	\begin{split}
		e(t) 
		:= &
		Y_{t} - \bar{ Y} _{ t }  \\
		= &
		e^{ - \lambda ( t - t_{n} ) }
		e ( t_{ n } )
		+
		\int_{ t_{ n } }^{ t }
		e^{ - \lambda ( t - r ) }
		F ( Y_{ r } ) 
		- 
		e^{ - \lambda ( t - t_{n } ) }
		F^{ ( \tau ) } ( \bar{Y}_{ t_{ n } } ) 
		\dd r \\
		&
		+
		\int_{ t_{ n } }^{ t }
		\big( e^{ - \lambda ( t - r ) }
		-
		e^{ - \lambda ( t - t_{ n } ) }	\big)
		b ( 1 - \theta ) \dd W_{ r }.
	\end{split}
\end{equation*}
It is straightforward to check that $ e(s) $
is the solution to the SDE
\begin{equation*}
	\dd e(t) =
	\big[ -\lambda e(t) + F( Y_{t} )
	- e^{- \lambda ( t - t_{n} ) }
	F^{ ( \tau ) } ( \bar{Y}_{ t_{n} }  )
	\big]  \dd t
	+
	\big( 1 - e^{ - \lambda( t - t_{n} )  }
	\big) b ( 1 - \theta )  \dd W_{ t }.
\end{equation*}
Using It\^o's formula, we have  
\begin{equation}\label{eq:tmp}
	\begin{split}
		| e( t ) |^2
		= &
		| e( t_{ n } ) |^2
		+
		\int_{ t_{ n } }^{ t }
		\big| ( 1 - e^{ - \lambda ( r - t_{ n } ) } )
		b ( 1 - \theta ) \big|^2 \dd r  \\
		& +
		2 \int_{ t_{ n } }^{ t }
		\left\langle e( r ),
		- \lambda e(r) + F( Y_{r} )
		- e^{- \lambda ( r - t_{n} ) }
		F^{ ( \tau ) } ( \bar{Y}_{ t_{n}  } )
		\right\rangle  \dd r \\
		& +
		2 \int_{ t_{ n } }^{ t }
		\left\langle e( r ),
		( 1 - e^{ - \lambda( r - t_{n} )  }
		) b ( 1 - \theta )  
		\right\rangle  \dd W_{ r }.
	\end{split}
\end{equation}
By Taylor's Theorem, there exists
a constant $ K $ such that
$ 1 - e^{ - \lambda ( r - t_{ n } ) } 
\leq  \lambda ( r - t_{ n } ) + K ( r - t_{ n } )^2 $.
Taking the expectation of both sides of \eqref{eq:tmp} and using It\^o's Isometry
\begin{equation}\label{eq:one-step-error}
	\begin{split}
		\E [ | e( t ) |^2 ] 
		= &
		\E [ | e( t_{ n } ) |^2 ]
		+
		2 \int_{ t_{ n } }^{ t }
		\E \big[ 
		\langle e( r ),
		- \lambda e(r) 	\rangle \big]
		+ 
		\E \big[ 
		\langle e( r ),
		F( Y_{r} ) 
		- e^{- \lambda ( r - t_{n} ) }
		F^{ ( \tau ) } ( \bar{Y}_{ t_{n} }  )
		\rangle \big]  \dd r \\
		& +
		\int_{ t_{ n } }^{ t }
		\E \big[ 
		\big| ( 1 - e^{ - \lambda ( r - t_{ n } ) } )
		b ( 1 - \theta ) \big|^2  \big] \dd r  \\
		\leq
		&
		\E [ | e( t_{ n } ) |^2 ]
		+
		2 \int_{ t_{ n } }^{ t }
		\E \big[ 
		\langle e( r ),
		F( Y_{r} ) 
		- F( \bar{Y}_{r} ) 
		\rangle \big]  \dd r  \\
		& +
		2 \int_{ t_{ n } }^{ t }
		\E \big[ 
		\langle e( r ),
		F( \bar{Y}_{r} ) 
		- e^{- \lambda ( r - t_{n} ) }
		F^{ ( \tau ) } ( \bar{Y}_{ t_{n} }  )
		\rangle \big]  \dd r
		+
		C
		\int_{ t_{ n } }^{ t }
		( r - t_{ n } )^2
		+ ( r - t_{ n } )^4  \dd r  \\
		\leq  
		&
		\E [ | e( t_{ n } ) |^2 ]
		+
		2 J_1 + 2 J_2 
		+ C \tau^3.
	\end{split}
\end{equation}
For $ J_1 $, according to the definition
of function $ F $, we have
\begin{equation*}
	J_1 \leq C 
	\int_{ t_{ n } }^{ t }
	\E \big[ | e( r ) |^2
	\big]  \dd r.
\end{equation*}
For $ J_2 $,
\begin{equation*}
	\begin{split}
		J_2 
		= &
		\int_{ t_{ n } }^{ t }
		\E \big[ 
		\langle e( r ),
		F ( \bar{Y}_{r} )
		- 
		F ( \bar{Y}_{ t_{ n } } )
		\rangle \big]  \dd r
		+
		\int_{ t_{ n } }^{ t }
		\E \big[ 
		\langle e( r ),
		F ( \bar{Y}_{ t_{ n } } )
		- e^{- \lambda ( r - t_{n} ) }
		F^{ ( \tau ) } ( \bar{Y}_{ t_{n} }  )
		\rangle \big]  \dd r \\
		= & :
		J_{ 21 } + J_{ 22}.
	\end{split}
\end{equation*}
For $ J_{ 22 } $, using the Young inequality,
Lemma \ref{lem:bound-moment-numerical}
and Lemma \ref{lem:inverse-bound-moment-numerical},
we obtain
\begin{equation*}
	\begin{split}
		J_{ 22 }
		= &
		\int_{ t_{ n } }^{ t }
		\E \big[ 
		\langle e( r ),
		( 1 - e^{- \lambda ( r - t_{n} ) } )
		F ( \bar{Y}_{ t_{ n } } )
		\rangle \big]  \dd r
		+
		\int_{ t_{ n } }^{ t }
		\E \big[ 
		\langle e( r ),
		e^{- \lambda ( r - t_{n} ) }
		\big( F ( \bar{Y}_{ t_{ n } } )
		- 
		F^{ ( \tau ) } ( \bar{Y}_{ t_{n} }  ) \big)
		\rangle \big]  \dd r \\
		\leq 
		&
		\int_{ t_{ n } }^{ t }
		\E \big[ | e( r ) |^2 \big] \dd r
		+ C \tau^3
		\E \big[ | 	F ( \bar{Y}_{ t_{ n } } ) |^2 
		\big] 
		+
		C \tau^3 
		\E \Big[ \Big| 	
		\frac{ \tau | F ( \bar{Y}_{ t_{ n } } )  |^2
			F ( \bar{Y}_{ t_{ n } } ) }
		{ 1 + \tau | F ( \bar{Y}_{ t_{ n } } ) |^2 } 
		\Big|^2  \Big]   \\
		\leq 
		& 
		 \int_{ t_{ n } }^{ t }
		\E \big[ | e( r ) |^2 \big] \dd r
		+
		C \tau^3 \big( 1 + 
		\E \big[ | \bar{Y}_{ t_{ n } } |^{ 
			\frac{2 ( \gamma - \theta)}{1 - \theta } }
		]
		+
		\E \big[ 
		| \bar{Y}_{ t_{ n } } |^{ 
			\frac{2 ( \theta + 1 )}{ \theta - 1  } }
		\big]
		+
		\E \big[ 
		| \bar{Y}_{ t_{ n } } |^{ 
			\frac{6 ( \gamma - \theta)}{1 - \theta } }
		\big]
		+
		\E \big[ 
		| \bar{Y}_{ t_{ n } } |^{ 
			\frac{6 ( \theta + 1 )}{ \theta - 1  } }
		\big]
		\big) \\
		\leq 
		&
		\int_{ t_{ n } }^{ t }
		\E \big[ | e( r ) |^2 \big] \dd r
		+ 
		C \tau^3.
	\end{split}
\end{equation*}
Using Taylor theorem, we get
\begin{equation*}
	\begin{split}
		& \E \big[ 
		\langle e( r ),
		F ( \bar{Y}( r ) )
		- 
		F ( \bar{Y}_{ t_{ n } } )
		\rangle \big]   \\
		= &
		\E \big[ 
		\langle e( r ),
		F' ( \bar{Y}_{ t_{ n } } )
		\int_{ t_{ n } }^{ r }
		e^{ - \lambda ( s - t_{ n } ) }
		F^{ ( \tau ) } ( \bar{Y}_{ t_{n} }  )
		\dd s
		\rangle \big]  
		+
		\E \big[ 
		\langle e( r ),
		F' ( \bar{Y}_{ t_{ n } } )
		\int_{ t_{ n } }^{ r }
		e^{ - \lambda ( s - t_{ n } ) }
		b ( 1 - \theta ) \dd W_s
		\rangle \big]  	\\
		& +
		\E \big[ 
		\langle e( r ),
		F' ( \bar{Y}_{ t_{ n } } )
		(1 - e^{ - \lambda ( r - t_{ n } ) } )
		 \bar{Y}_{ t_{ n } }
		\rangle \big] 
		+
		\E \big[ 
		\langle e( r ),
		\int_{0}^{1}
		F'' ( \bar{Y}_{ t_{ n } } 
		+ \sigma ( \bar{Y}_{ r } 
		- \bar{Y}_{ t_{ n } } )  )
		| \bar{Y}_{ r } 
		- \bar{Y}_{ t_{ n } } |^2 \dd \sigma
		\rangle \big] \\
		= & :
		H_1 + H_2 + H_3 + H_4.
	\end{split}
\end{equation*}
Note that 
\begin{equation*}
	\begin{split}
		\bar{Y}_{ r } 
		- \bar{Y}_{ t_{ n } } 
		= &
		\bar{Y}_{ r } 
		-
		e^{ - \lambda ( r - t_{ n } ) }
		\bar{Y}_{ t_{ n } } 
		- ( 1 -  e^{ - \lambda ( r - t_{ n } ) } )
		\bar{Y}_{ t_{ n } } \\
		= &
		\int_{ t_{ n } }^{ r }
		e^{ - \lambda ( s - t_{ n } ) }	
		F^{ ( \tau ) } ( \bar{Y}_{ t_{ n } } ) \dd s 
		+
		\int_{ t_{ n } }^{ r }
		e^{ - \lambda ( s - t_{ n } ) }	
		b ( 1 - \theta ) \dd W_{s} 
		- ( 1 -  e^{ - \lambda ( r - t_{ n } ) } )
		\bar{Y}_{ t_{ n } } .	
	\end{split}
\end{equation*}
For $ p \geq 2 $,
it's straightforward to derive that
\begin{equation}\label{eq:numerical-errpr}
	\begin{split}
		& \E [ | \bar{Y}_{ r } 
		- \bar{Y}_{ t_{ n } } |^p  ] \\
		\leq &
		C_{ p } \E \Big[ 
		\Big| \int_{ t_{ n } }^{ r }
		e^{ - \lambda ( s - t_{ n } ) }	
		F^{ ( \tau ) } ( \bar{Y}_{ t_{ n } } ) \dd s 
		\Big|^p \Big] 
		+
		C_{ p } \E \Big[ \Big|
		\int_{ t_{ n } }^{ r }
		e^{ - \lambda ( s - t_{ n } ) }	
		b ( 1 - \theta ) \dd W_{s} \Big|^p
		\Big]
		+ C_{ p } \E [
		| ( 1 -  e^{ - \lambda ( r - t_{ n } ) } )
		\bar{Y}_{ t_{ n } } |^p ] \\
		\leq 
		&
		C_{ p } \tau^{ p - 1 }
		 \int_{ t_{ n } }^{ r }
		\E [ | \bar{Y}_{ t_{ n } } |^p ] \dd s
		+ 
		C_{ p }  \Big( 
		\int_{ t_{ n } }^{ r }	
		| b ( 1 - \theta ) |^2  \dd s 
		\Big)^{\frac{ p }{ 2 } } 
		+ C_{ p } \tau^{ p }
		\E [
		| 	\bar{Y}_{ t_{ n } } |^p ] \\
		\leq 
		&
		C \tau^2 
		\big( 1 + \E \big[ | \bar{Y}_{ t_{ n } } |^{ 
			\frac{2 ( \gamma - \theta)}{1 - \theta } }
		]
		+
		\E \big[ 
		| \bar{Y}_{ t_{ n } } |^{ 
			\frac{2 ( \theta + 1 )}{ \theta - 1  } }
		\big] \big)
		+ 
		C \tau^{ \frac{ p }{ 2 } } \\
		\leq &
		C \tau^{ \frac{ p }{ 2 } }.
	\end{split}
\end{equation}
Using the Young inequality and
Lemma \ref{lem:bound-moment-numerical},
\ref{lem:inverse-bound-moment-numerical},
we obtain 
\begin{equation*}
	\begin{split}
	H_1 
	\leq &
	\E [ | e( r ) |^2 ]
	+
	C 
	\E 
	\Big[ 
	| F' ( \bar{Y}_{ t_{ n } } ) |^2
	\Big| \int_{ t_{ n } }^{ r }
	e^{ - \lambda ( s - t_{ n } ) }
	F^{ ( \tau ) } ( \bar{Y}_{ t_{n} }  )
	\dd s \Big|^2
	\Big]  \\
	\leq
	&
	\E [ | e( r ) |^2 ]
	+
	C \tau
	\E  \big[ 
	| F' ( \bar{Y}_{ t_{ n } } ) |^2
      \int_{ t_{ n } }^{ r }
	 | F^{ ( \tau ) } ( \bar{Y}_{ t_{n} }  ) |^2
	\dd s 
	\big]  \\
	\leq
	&
	\E [ | e( r ) |^2 ]
	+
	C \tau^2
	\big( \E [ | F' ( \bar{Y}_{ t_{ n } } ) |^4 ] 
	\big)^{ \frac{ 1 } { 2 } }
		\big( \E [ | F ( \bar{Y}_{ t_{ n } } ) |^4 ] 
	\big)^{ \frac{ 1 } { 2 } } \\
	\leq 
	&
	\E [ | e( r ) |^2 ]
	+
	C \tau^2 
	\big( 1 + \E \big[ | \bar{Y}_{ t_{ n } } |^{ 
		\frac{4 ( \gamma - \theta)}{1 - \theta } }
	]
	+
	\E \big[ 
	| \bar{Y}_{ t_{ n } } |^{ 
		\frac{4 ( \theta + 1 )}{ \theta - 1  } }
	\big] \big)^{ \frac{ 1 } { 2 } }
	\big( 1 + \E \big[ | \bar{Y}_{ t_{ n } } |^{ 
		\frac{4 ( \gamma - 1)}{1 - \theta } }
	]
	+
	\E \big[ 
	| \bar{Y}_{ t_{ n } } |^{ 
		\frac{4 ( 2 )}{ \theta - 1  } }
	\big] \big)^{ \frac{ 1 } { 2 } } \\
	\leq
	&
	\E [ | e( r ) |^2 ]
	+
	C \tau^2 .
	\end{split}
\end{equation*}
For $ H_2 $, by \eqref{eq:tmp}
and
\eqref{eq:numerical-errpr},
we have
\begin{equation*}
	\begin{split}
		H_2 
		= &
		\E \big[ 
		\langle 
		e^{ - \lambda ( r - t_{ n } ) } e( t_{ n } ),
		F' ( \bar{Y}_{ t_{ n } } )
		\int_{ t_{ n } }^{ r }
		e^{ - \lambda ( s - t_{ n } ) }
		b ( 1 - \theta ) \dd W_s
		\rangle \big]  	\\
		& + 
		\E \Big[ 
		\Big\langle 
		\int_{t_{n} }^{ r }
		e^{ - \lambda ( r - \sigma )}
		F( Y_{ \sigma } ) 
		-
		e^{ - \lambda ( r - t_{n} )}
		F^{ ( \tau ) }( \bar{Y}_{ t_{n} } )
		\dd \sigma ,
		F' ( \bar{Y}_{ t_{ n } } )
		\int_{ t_{ n } }^{ r }
		e^{ - \lambda ( s - t_{ n } ) }
		b ( 1 - \theta ) \dd W_s
		\Big\rangle \Big] \\
		& +
		\E \Big[ 
		\Big\langle 
		\int_{t_{n} }^{ r }
		\big( 
		e^{ - \lambda ( r - \sigma )}
		-
		e^{ - \lambda ( r - t_{n} )}
		\big)
		b ( 1 - \theta ) \dd W_{ \sigma },
		F' ( \bar{Y}_{ t_{ n } } )
		\int_{ t_{ n } }^{ r }
		e^{ - \lambda ( s - t_{ n } ) }
		b ( 1 - \theta ) \dd W_s
		\Big\rangle \Big] \\
		= &
		\E \Big[ 
		\Big\langle 
		\int_{t_{n} }^{ r }
		( e^{ - \lambda ( r - \sigma )}
		- e^{ - \lambda ( r - t_{n} )} )
		F( Y_{ \sigma } ) 
		\dd \sigma ,
		F' ( \bar{Y}_{ t_{ n } } )
		\int_{ t_{ n } }^{ r }
		e^{ - \lambda ( s - t_{ n } ) }
		b ( 1 - \theta ) \dd W_s
		\Big\rangle \Big] \\
		& +
		\E \Big[ 
		\Big\langle 
		\int_{t_{n} }^{ r }
		e^{ - \lambda ( r - t_{n} )}
		\big( F( Y_{ \sigma } ) 
		-
		F( Y_{ t_{ n } } ) \big)
		\dd \sigma ,
		F' ( \bar{Y}_{ t_{ n } } )
		\int_{ t_{ n } }^{ r }
		e^{ - \lambda ( s - t_{ n } ) }
		b ( 1 - \theta ) \dd W_s
		\Big\rangle \Big] \\
		& +
		\E \Big[ 
		\Big\langle 
		\int_{t_{n} }^{ r }
		e^{ - \lambda ( r - t_{ n } )}
		\big( 
		F( Y_{ t_{ n } } ) 
		-
		F( \bar{Y}_{ t_{n} } ) \big)
		\dd \sigma ,
		F' ( \bar{Y}_{ t_{ n } } )
		\int_{ t_{ n } }^{ r }
		e^{ - \lambda ( s - t_{ n } ) }
		b ( 1 - \theta ) \dd W_s
		\Big\rangle \Big] \\
		& +
		\E \Big[ 
		\Big\langle 
		\int_{t_{n} }^{ r }
		e^{ - \lambda ( r - t_{ n } )}
		\big( F( \bar{Y}_{ t_{n } } ) 
		-
		F^{ ( \tau ) }( \bar{Y}_{ t_{n} } ) 
		\big)	\dd \sigma ,
		F' ( \bar{Y}_{ t_{ n } } )
		\int_{ t_{ n } }^{ r }
		e^{ - \lambda ( s - t_{ n } ) }
		b ( 1 - \theta ) \dd W_s
		\Big\rangle \Big] \\
		& +
		\E \Big[ 
		\Big\langle 
		\int_{t_{n} }^{ r }
		\big( 
		e^{ - \lambda ( r - \sigma )}
		-
		e^{ - \lambda ( r - t_{n} )}
		\big)
		b ( 1 - \theta ) \dd W_{ \sigma },
		F' ( \bar{Y}_{ t_{ n } } )
		\int_{ t_{ n } }^{ r }
		e^{ - \lambda ( s - t_{ n } ) }
		b ( 1 - \theta ) \dd W_s
		\Big\rangle \Big].
	\end{split}
\end{equation*}
Then, using H\"older inequality,
the It\^o isometry and the bounded moment
of exact and numerial solution to,
one can get \eqref{eq:Lamperti}
\begin{equation*}
	\begin{split}
	H_2 \leq &
	\Big( \E \Big[ 
	\Big| 
	\int_{t_{n} }^{ r }
	( e^{ - \lambda ( r - \sigma )}
	- e^{ - \lambda ( r - t_{n} )} )
	F( Y_{ \sigma } ) 
	\dd \sigma \Big|^2 \Big]  
	\E \Big[ \Big| 
	F' ( \bar{Y}_{ t_{ n } } )
	\int_{ t_{ n } }^{ r }
	e^{ - \lambda ( s - t_{ n } ) }
	b ( 1 - \theta ) \dd W_s
	\Big|^2 \Big] \Big)^{\frac{1}{2}} \\
	& +
	\Big( \E \Big[ 
	\Big| 
	\int_{t_{n} }^{ r }
	e^{ - \lambda ( r - t_{n} )}
	( F( Y_{ \sigma } ) 
	-
	F( Y_{ t_{n} } ) )
	\dd \sigma \Big|^2 \Big]  
	\E \Big[ \Big| 
	F' ( \bar{Y}_{ t_{ n } } )
	\int_{ t_{ n } }^{ r }
	e^{ - \lambda ( s - t_{ n } ) }
	b ( 1 - \theta ) \dd W_s
	\Big|^2 \Big] \Big)^{\frac{1}{2}} \\
	& +
	\Big( \E \Big[ 
	\Big| 
	\int_{t_{n} }^{ r }
	( e^{ - \lambda ( r - \sigma )}
	- e^{ - \lambda ( r - t_{n} )} )
	b ( 1 - \theta )
	\dd W_{\sigma} \Big|^2 \Big]  
	\E \Big[ \Big| 
	F' ( \bar{Y}_{ t_{ n } } )
	\int_{ t_{ n } }^{ r }
	e^{ - \lambda ( s - t_{ n } ) }
	b ( 1 - \theta ) \dd W_s
	\Big|^2 \Big] \Big)^{\frac{1}{2}} \\
	\leq &
	C \tau^{ \frac{ 3 }{ 2 } }
	\Big( \E \Big[  
	\int_{t_{n} }^{ r }
	|
	F( Y_{ \sigma } )  |^2
	\dd \sigma  \Big]  
	\cdot
	\E \Big[ \big| 
	F' ( \bar{Y}_{ t_{ n } } ) \big|^2
	\int_{ t_{ n } }^{ r }
	|
	b ( 1 - \theta ) |^2 \dd s
	\Big] \Big)^{\frac{1}{2}} \\
	& +
	C \tau^{ \frac{1}{2}}
	\Big( \E \Big[ 
	\int_{t_{n} }^{ r }
	| ( F( Y_{ \sigma } ) 
	-
	F( Y_{ t_{n} } ) ) |^2
	\dd \sigma  \Big]  
	\cdot
	\E \Big[ \big| 
	F' ( \bar{Y}_{ t_{ n } } ) \big|^2
	\int_{ t_{ n } }^{ r }
	|
	b ( 1 - \theta ) |^2  \dd s
	\Big] \Big)^{\frac{1}{2}} \\
	& +
	C \Big( \E \Big[  
	\int_{t_{n} }^{ r }
	\big| ( e^{ - \lambda ( r - \sigma )}
	- e^{ - \lambda ( r - t_{n} )} ) \big|^2
	\dd \sigma \Big|^2 \Big]  
	\E \Big[ \big| 
	F' ( \bar{Y}_{ t_{ n } } ) \big|^2
	\int_{ t_{ n } }^{ r }
	| b ( 1 - \theta ) |^2 \dd s
	\Big] \Big)^{\frac{1}{2}} \\
	\leq &
	C \tau^{ \frac{ 5 }{ 2 } }
	\Big( \E \Big[  
	\sup_{ t_{ n } \leq \sigma \leq t_{ n+1 }}
	|
	F( Y_{ \sigma } )  |^2  \Big]  
	\cdot
	\E \Big[ \big| 
	F' ( \bar{Y}_{ t_{ n } } ) \big|^2
	\Big] \Big)^{\frac{1}{2}} \\
	& +
	C \tau^{ \frac{1}{2}}
	\Big( \E \Big[ 
	\int_{t_{n} }^{ r }
	| F' ( \xi_{ n } ) |^2
	| Y_{ \sigma } -
	Y_{ t_{n} }  |^2
	\dd \sigma  \Big]  
	\cdot
	\E \Big[ \big| 
	F' ( \bar{Y}_{ t_{ n } } ) \big|^2
	\Big] \Big)^{\frac{1}{2}} 
	+
	C \tau^2  \Big( \E \Big[  
	\big| 
	F' ( \bar{Y}_{ t_{ n } } ) \big|^2
	\Big] \Big)^{\frac{1}{2}} \\
	\leq 
	&
	C \tau^{\frac{5}{2}}
	+
	C\tau^2,	
	\end{split}
\end{equation*}
where $ \xi_{ n } = Y_{ t_{ n } } 
 + \xi \big( Y_{\sigma} - Y_{ t_{ n } } \big)$,
 $ \xi \in ( 0, 1) $.
 
For $ H_3 $, the Cauchy-Schwarz 
inequality implies
\begin{equation*}
	\begin{split}
	H_3 
	\leq
	&
	\E \big[ 
	| e( r )|^2 \big] 
	+
	\E \big[ 
	| F' ( \bar{Y}_{ t_{ n } } )
	(1 - e^{ - \lambda ( s - t_{ n } ) } )
	b ( 1 - \theta ) \bar{Y}_{ t_{ n } } |^2
	\big]  \\
	\leq 
	&
	\E \big[ 
	| e( r )|^2 \big] 
	+ C \tau^2 
	\big( \E [ | F' ( \bar{Y}_{ t_{ n } } ) |^4 ]
	\big)^{ \frac{ 1 }{ 2 } } 
	\big( \E [ | \bar{Y}_{ t_{ n } } |^4 ]
	\big)^{ \frac{ 1 }{ 2 } } \\
	\leq 
	&
	\E \big[ 
	| e( r )|^2 \big] 
	+ C \tau^2 .
	\end{split}
\end{equation*}
For $ H_4 $, using the the Cauchy-Schwarz 
inequality
and \eqref{eq:numerical-errpr}, we have
\begin{equation*}
	\begin{split}
		H_4
		& \leq
		\E \big[ 
		| e( r )|^2 \big] 
		+
		\E \big[ 
		| F'' \big( \bar{Y}_{ t_{ n } }
		+ \sigma ( \bar{Y}_{r} -\bar{Y}_{ t_{ n } } ) \big)
		| \bar{Y}_{r} -\bar{Y}_{ t_{ n } } |^2  |^2
		\big] \\
		& \leq 
		\E \big[ 
		| e( r )|^2 \big] 
		+ C  \Big( \E [ |  F'' \big( \bar{Y}_{ t_{ n } }
		+ \sigma ( \bar{Y}_{r} -\bar{Y}_{ t_{ n } } ) \big) |^4]
		\Big)^{ \frac{1}{2} }
		\big( \E [ | \bar{Y}_{r} -\bar{Y}_{ t_{ n } } |^8]
		\big)^{\frac{1}{2}} \\
		& \leq
		\E \big[ 
		| e( r )|^2 \big] 
		+C \tau^2.
	\end{split}
\end{equation*}
Note that,
\begin{equation*}
		J_{ 21 }
		=
		\int_{ t_{ n }}^{t}
		\E \big[ 
		\langle e( r ),
		F ( \bar{Y}( r ) )
		- 
		F ( \bar{Y}_{ t_{ n } } )
		\rangle \big] \dd r 
		= 
		\int_{ t_{ n }}^{t}
		H_1 + H_2 + H_3 + H_4 \dd r
		\leq 
		3 \int_{ t_{ n }}^{t}
	     \E [ | e( r ) |^2 ] \dd r
		+ C \tau^3.
\end{equation*}
Therefore, from above discuss, we have
$$
J_2 \leq
 4 \int_{t_n }^{ t } 
\E [ | e( r ) |^2 ] \dd r 
+
C \tau^3.
$$
Substituting $ J_1 $
 and $ J_2 $ 
into equation \eqref{eq:one-step-error},
 we can get
\begin{equation*}
	\begin{split}
		\E [ | e( t ) |^2 ]
		\leq
		&
		\E [ | e( t_{ n } ) |^2 ]
		+
		C \int_{ t_{ n } }^{ t }
		\E \big[ 
		| e( r ) |^2 \big]  \dd r
		+
		C \tau^3 .
	\end{split}
\end{equation*}
%
Thus, for any $t\in[0,T]$ we have that 
\begin{equation*}
	\begin{split}
		\E [ | e( t ) |^2 ]
		\leq
		&
		C \int_{ 0 }^{ t }
		\E \big[ 
		| e( r ) |^2 \big]  \dd r
		+
		C \tau^2.
	\end{split}
\end{equation*}
Use Gronwall's inequality, we get
\begin{equation*}
	\begin{split}
		\E [ | e( t ) |^2 ]
		\leq
		&
		C \tau^2.
	\end{split}
\end{equation*}
This completes the proof.
\qed

\section{Probability of positivity}
\label{sec:positivity}
We now examine, following closely the analysis in \cite{KLM2020}, the probability of solutions of 
\eqref{eq:numerical-backstop} 
becoming negative after a single 
step with this strategy, 
and hence triggering a use of the backstop method.  
This probability is given by
$$
\mathbb{P}\left[Y_{k+1}<0 \mid 
Y_{k}=y>0\right]= \Phi(a(y)),
$$
where 
$$
\Phi(x)=\frac{1}{\sqrt{2 \pi}} 
\int_{-\infty}^{x} e^{-s^{2} / 2} \dd s
\quad \text{and} \quad 
a(y) := \frac{ y+ 
	 \Ftau(y) }
	 { b ( \theta - 1 ) \sqrt{ \tau } }.
$$

\begin{theorem}\label{thm:probability-positive} 
	Let $\left\{Y_{n}\right\}_{n=0}^{N}$ 
	be a solution of 
	\eqref{eq:numerical-backstop}, 
	with initial value 
	$Y_{0}>0$. 
	Suppose also that 
	$Y_{0} \leq M_1 $. 
	Then, for each $\varepsilon \in(0,1)$ 
	there exists $ \tau >0 $ such that
$$
\mathbb{P}\left[\mathcal{R}_{N}\right]>1-\varepsilon
$$
where 
$\mathcal{R}_{N}:=\bigcap_{j=0}^{N}\left\{Y_{j}>0\right\}$.
\end{theorem}

\textbf{Proof of Theorem \ref{thm:probability-positive}}
Since, 
by \eqref{eq:numerical-backstop}, 
the backstop method will ensure 
positivity over a single step 
$ \tau $, 
the event $ \left\{Y_{n+1}>0\right\}$ 
is equivalent to the following:
$$
\begin{aligned}
   \bigg\{ \frac{\Delta W_{n}}{\sqrt{ \tau }}
	<\frac{1}{b ( \theta - 1 ) }
	\left( \frac{Y_{n}}{\sqrt{\tau}}
	+\frac{ F(Y_{ n }) \sqrt{ \tau } }
	{ 1 + \tau | F( Y_{ n }) |^2 } \right)
	\bigg\}.
\end{aligned}
$$
From Lemmas \ref{lem:bound-moment-numerical},
\ref{lem:inverse-bound-moment-numerical},
we can easy to check that there exist
constants 
$ M_1 $, $ M_2 $ such that
\begin{equation*}
	M_1 \leq 
	\sup_{n \in \mathbb{N}} | Y_n | 
	\leq M_2,
	\quad  a.s.
\end{equation*}
Thus, we have
\begin{equation*}\label{eq:esti-I}
	\begin{split}
		\Ftau \sqrt{ \tau }
		  = 
		  &
		\frac{ 
			F(Y_{ n } )  \sqrt{ \tau } }
		{ 1 + \tau | F( Y_{ n } ) |^2} \\
		= &
		\frac{
			( \theta - 1 ) \sqrt{ \tau }
			\big[ a_2 Y_{ n }^{ 
				\frac{\gamma -  \theta }
				{ 1 - \theta } } - 
			a_{ - 1 } Y_{  n }^{ 
				\frac{  \theta + 1  }{ \theta - 1  } }
			+ a_{ 0 } 
			Y_{ n }^{ 
				\frac{  \theta  }
				{ \theta - 1 } }
			+ \frac{ 1 }{ 2 }
			b^2 \theta Y_{ n }^{ -1 }
			\big] }
		{ 1 + \tau ( \theta - 1 )^2 \Big| 
			a_2 Y_{ n }^{ 
				\frac{\gamma  -  \theta }
				{ 1 - \theta } } - 
			a_{ - 1 } Y_{ n }^{ 
				\frac{  \theta + 1 }{ \theta - 1  } }
			+ a_{ 0 } 
			Y_{ n }^{ 
				\frac{  \theta  }
				{ \theta - 1 } }
			+ \frac{ 1 }{ 2 }
			b^2 \theta Y_{ n }^{-1}
			\Big|^2 }		\\
		\geq  &
		\frac{
			( \theta - 1 ) \sqrt{ \tau }
			\big[ a_2 M_{ 2 }^{ 
				\frac{\gamma -  \theta }
				{ 1 - \theta } } - 
			a_{ - 1 } M_{ 2 }^{ 
				\frac{  \theta + 1  }{ \theta - 1  } }
			+ a_{ 0 } 
			M_{ 1 }^{ 
				\frac{  \theta  }
				{ \theta - 1 } }
			+ \frac{ 1 }{ 2 }
			b^2 \theta M_{ 2 }^{ -1 }
			\big] }
		{ 1 + 4 \tau ( \theta - 1 )^2 \Big(
			a_2^2 M_{ 1 }^{ 
				\frac{2\gamma  -  2\theta }
				{ 1 - \theta } } + 
			a_{ - 1 }^2 M_{ 2 }^{ 
				\frac{  2\theta + 2 }{ \theta - 1  } }
			+ a_{ 0 }^2 
			M_{ 2 }^{ 
				\frac{  2\theta  }
				{ \theta - 1 } }
			+ \frac{ 1 }{ 4 }
			b^4 \theta M_{ 1 }^{-2}
			 \Big)
		}		.
	\end{split}
\end{equation*}
To prove 
$ Y_{ n + 1 } > 0 $, 
it is sufficient to prove
\begin{equation*}
	\begin{split}
		\frac{\Delta W_{n}}
		{\sqrt{ \tau }}
		<
		\frac{1}{b ( \theta - 1 ) }
		\left( \frac{Y_{n}}{\sqrt{\tau}}
		+
	\frac{
		( \theta - 1 ) \sqrt{ \tau }
		\big[ a_2 M_{ 2 }^{ 
			\frac{\gamma -  \theta }
			{ 1 - \theta } } - 
		a_{ - 1 } M_{ 2 }^{ 
			\frac{  \theta + 1  }{ \theta - 1  } }
		+ a_{ 0 } 
		M_{ 1 }^{ 
			\frac{  \theta  }
			{ \theta - 1 } }
		+ \frac{ 1 }{ 2 }
		b^2 \theta M_{ 2 }^{ -1 }
		\big] }
	{ 1 + 4 \tau ( \theta - 1 )^2 \Big(
		a_2^2 M_{ 1 }^{ 
			\frac{2\gamma  -  2\theta }
			{ 1 - \theta } } + 
		a_{ - 1 }^2 M_{ 2 }^{ 
			\frac{  2\theta + 2 }{ \theta - 1  } }
		+ a_{ 0 }^2 
		M_{ 2 }^{ 
			\frac{  2\theta  }
			{ \theta - 1 } }
		+ \frac{ 1 }{ 4 }
		b^4 \theta M_{ 1 }^{-2}
		\Big)
	} \right).
	\end{split}
\end{equation*}
Let $ \Phi $ denote the
distribution function of a 
standard Normal random variable, 
and suppose 
$ \left\{ \xi_{ n } \right\}_{n \in \mathbb{ N } } $ 
is a sequence of mutually 
independent standard normal random variables. 
Then
{\small
\begin{equation*}
	\begin{split}
	& \mathbb{P}
	\left[ Y_{n+1} > 0 
	\mid \mathcal{R}_{n} \right] \\
	& =
	\frac{ \mathbb{P}\left[
		\left\{ Y_{n+1} > 0 \right\} \cap 
		\left\{\mathcal{R}_{n} \right\}\right]}
	{\mathbb{ P }
		\left[\mathcal{R}_{n} \right]}
	=
	\frac{\mathbb{E}\left[\mathbb{E}
		\left[I_{\left\{Y_{n+1}>0\right\} 
			\cap\left\{\mathcal{R}_{n} \right\}} \mid \mathcal{F}_{t_{n}}\right]\right]}
	{\mathbb{P}\left[\mathcal{R}_{n} \right]}\\
	&
	\geq 
	\mathbb{E}\left[
	\mathbb{E} \left[ I_{\left\{ Y_{n+1} > 0 \right\} 
		\cap
		\left\{\mathcal{R}_{n} \right\}} 
	\mid \mathcal{F}_{t_{n}}\right]
	\right]
	=
	\mathbb{E}\left[\mathbb{P}
	\left[\left\{Y_{n+1}>0\right\} 
	\cap
	\left\{\mathcal{R}_{n} 
	\right\} \mid \mathcal{F}_{t_{n}}\right]\right]\\
	&
	=
	\mathbb{E}\left[ \mathbb{P}
	\left[	\frac{\Delta W_{n}}
	{\sqrt{ \tau }}
	<
	\frac{1}{b ( \theta - 1 ) }
	\left( \frac{Y_{n}}{\sqrt{\tau}}
	+
	\frac{
		( \theta - 1 ) \sqrt{ \tau }
		\big[ a_2 M_{ 2 }^{ 
			\frac{\gamma -  \theta }
			{ 1 - \theta } } - 
		a_{ - 1 } M_{ 2 }^{ 
			\frac{  \theta + 1  }{ \theta - 1  } }
		+ a_{ 0 } 
		M_{ 1 }^{ 
			\frac{  \theta  }
			{ \theta - 1 } }
		+ \frac{ 1 }{ 2 }
		b^2 \theta M_{ 2 }^{ -1 }
		\big] }
	{ 1 + 4 \tau ( \theta - 1 )^2 \Big(
		a_2^2 M_{ 1 }^{ 
			\frac{2\gamma  -  2\theta }
			{ 1 - \theta } } + 
		a_{ - 1 }^2 M_{ 2 }^{ 
			\frac{  2\theta + 2 }{ \theta - 1  } }
		+ a_{ 0 }^2 
		M_{ 2 }^{ 
			\frac{  2\theta  }
			{ \theta - 1 } }
		+ \frac{ 1 }{ 4 }
		b^4 \theta M_{ 1 }^{-2}
		\Big)
	} \right)
	\mid \mathcal{F}_{t_{n}}\right]\right]  \\
		&
	=
    \mathbb{P}
	\left[	\xi_{ n + 1 }
	<
	\frac{1}{b ( \theta - 1 ) }
	\left( \frac{Y_{n}}{\sqrt{\tau}}
	+
	\frac{
		( \theta - 1 ) \sqrt{ \tau }
		\big[ a_2 M_{ 2 }^{ 
			\frac{\gamma -  \theta }
			{ 1 - \theta } } - 
		a_{ - 1 } M_{ 2 }^{ 
			\frac{  \theta + 1  }{ \theta - 1  } }
		+ a_{ 0 } 
		M_{ 1 }^{ 
			\frac{  \theta  }
			{ \theta - 1 } }
		+ \frac{ 1 }{ 2 }
		b^2 \theta M_{ 2 }^{ -1 }
		\big] }
	{ 1 + 4 \tau ( \theta - 1 )^2 \Big(
		a_2^2 M_{ 1 }^{ 
			\frac{2\gamma  -  2\theta }
			{ 1 - \theta } } + 
		a_{ - 1 }^2 M_{ 2 }^{ 
			\frac{  2\theta + 2 }{ \theta - 1  } }
		+ a_{ 0 }^2 
		M_{ 2 }^{ 
			\frac{  2\theta  }
			{ \theta - 1 } }
		+ \frac{ 1 }{ 4 }
		b^4 \theta M_{ 1 }^{-2}
		\Big)
	} \right) \right]  \\
	&
	=
	\Phi \left(
	\frac{1}{b ( \theta - 1 ) }
	\left( \frac{Y_{n}}{\sqrt{\tau}}
	+
	\frac{
		( \theta - 1 ) \sqrt{ \tau }
		\big[ a_2 M_{ 2 }^{ 
			\frac{\gamma -  \theta }
			{ 1 - \theta } } - 
		a_{ - 1 } M_{ 2 }^{ 
			\frac{  \theta + 1  }{ \theta - 1  } }
		+ a_{ 0 } 
		M_{ 1 }^{ 
			\frac{  \theta  }
			{ \theta - 1 } }
		+ \frac{ 1 }{ 2 }
		b^2 \theta M_{ 2 }^{ -1 }
		\big] }
	{ 1 + 4 \tau ( \theta - 1 )^2 \Big(
		a_2^2 M_{ 1 }^{ 
			\frac{2\gamma  -  2\theta }
			{ 1 - \theta } } + 
		a_{ - 1 }^2 M_{ 2 }^{ 
			\frac{  2\theta + 2 }{ \theta - 1  } }
		+ a_{ 0 }^2 
		M_{ 2 }^{ 
			\frac{  2\theta  }
			{ \theta - 1 } }
		+ \frac{ 1 }{ 4 }
		b^4 \theta M_{ 1 }^{-2}
		\Big)
	} \right) \right).
	\end{split}
\end{equation*}
}
Since 
$ Y_{0} > 0 $  and $ \Phi $ takes values on 
$[0,1]$, 
we have
\begin{equation*}
	\begin{split}
			\mathbb{P}
			\left[\mathcal{R}_{N}\right]
			& =
			\mathbb{P}
			\left[\bigcap_{n=0}^{N} 
			\mathcal{R}_{n} \right] \\
		&
		=
		\prod_{n=1}^{N} 
		\mathbb{P}
		\left[\mathcal{R}_{n} \mid \mathcal{R}_{n-1}, \ldots, \mathcal{R}_{0} \right] 
	=
	\prod_{n=0}^{N-1}
	 \mathbb{P}\left[Y_{n+1}>0 
	 \mid \mathcal{R}_{n} \right] \\
	& 
	\geq 
	\Phi \left(
	\frac{1}{b ( \theta - 1 ) }
	\left( \frac{Y_{n}}{\sqrt{\tau}}
	+
	\frac{
		( \theta - 1 ) \sqrt{ \tau }
		\big[ a_2 M_{ 2 }^{ 
			\frac{\gamma -  \theta }
			{ 1 - \theta } } - 
		a_{ - 1 } M_{ 2 }^{ 
			\frac{  \theta + 1  }{ \theta - 1  } }
		+ a_{ 0 } 
		M_{ 1 }^{ 
			\frac{  \theta  }
			{ \theta - 1 } }
		+ \frac{ 1 }{ 2 }
		b^2 \theta M_{ 2 }^{ -1 }
		\big] }
	{ 1 + 4 \tau ( \theta - 1 )^2 \Big(
		a_2^2 M_{ 1 }^{ 
			\frac{2\gamma  -  2\theta }
			{ 1 - \theta } } + 
		a_{ - 1 }^2 M_{ 2 }^{ 
			\frac{  2\theta + 2 }{ \theta - 1  } }
		+ a_{ 0 }^2 
		M_{ 2 }^{ 
			\frac{  2\theta  }
			{ \theta - 1 } }
		+ \frac{ 1 }{ 4 }
		b^4 \theta M_{ 1 }^{-2}
		\Big)
	} \right) \right)^{N} .
	\end{split}
\end{equation*}
For convenience, we denote
$$
\mathbf{Y}_{ \tau }
:=
\frac{Y_{n}}{\sqrt{\tau}}
+
\frac{
	( \theta - 1 ) \sqrt{ \tau }
	\big[ a_2 M_{ 2 }^{ 
		\frac{\gamma -  \theta }
		{ 1 - \theta } } - 
	a_{ - 1 } M_{ 2 }^{ 
		\frac{  \theta + 1  }{ \theta - 1  } }
	+ a_{ 0 } 
	M_{ 1 }^{ 
		\frac{  \theta  }
		{ \theta - 1 } }
	+ \frac{ 1 }{ 2 }
	b^2 \theta M_{ 2 }^{ -1 }
	\big] }
{ 1 + 4 \tau ( \theta - 1 )^2 \Big(
	a_2^2 M_{ 1 }^{ 
		\frac{2\gamma  -  2\theta }
		{ 1 - \theta } } + 
	a_{ - 1 }^2 M_{ 2 }^{ 
		\frac{  2\theta + 2 }{ \theta - 1  } }
	+ a_{ 0 }^2 
	M_{ 2 }^{ 
		\frac{  2\theta  }
		{ \theta - 1 } }
	+ \frac{ 1 }{ 4 }
	b^4 \theta M_{ 1 }^{-2}
	\Big)
}.
 $$
Fix  $ \varepsilon \in (0,1) $,
 then for fix step size $ \tau $, 
 we have 
 \begin{equation}\label{eq:Phi-probability}
 	\Phi \left(
 	\frac{ \mathbf{Y}_{ \tau } }
 	{ b ( \theta - 1 ) }
 \right)^{ N } 
 	\geq 
 	1 - \varepsilon.
 \end{equation}
 To 
 \eqref{eq:Phi-probability}, 
 apply the following inequality 
 due to \cite{sasvari1999tight}
 $$
 \frac{1}{\sqrt{2 \pi}} \int_{-x}^{x} 
 e^{-s^{2} / 2} d s>\sqrt{1-e^{-x^{2} / 2}}, 
 \quad x \in \mathbb{R}^{+}
 $$
 along with the fact that 
 $ N =  T / \tau $, 
 leading us to seek 
 $ \tau $ 
 so that
 $$
 \left(\frac{1}{2}
 +\frac{1}{2} 
 \sqrt{1-\exp 
 	\left(-
 	\frac{  \mathbf{Y}_{ \tau } }
 	{2 b ( \theta - 1 )^{2}}
 	\right)}\right)^{ \frac{ T }{ \tau } } 
   \geq 
    1-\varepsilon.
 $$
Thus we derive the bound
 \begin{equation*}
 		\tau(\varepsilon)
 		:=
 		\sup \left\{ \tau \in (0,1): 
 	\mathbf{Y}_{ \tau } 
 	 \geq 
 	 \sqrt{\ln \left(1-\left(2(1-\varepsilon)^{\frac{h}{\rho T}}-1\right)^{2}\right)^{-2 b ( \theta - 1 )^{2}}} \right\} .	
 \end{equation*}
 $\tau (\varepsilon)$ 
 is uniquely defined for each 
 $\varepsilon \in(0,1)$ 
 because
\begin{equation*}
g( \tau )
:=
\mathbf{Y}_{ \tau }
-
\sqrt{-2 b ( \theta - 1 )^{2} 
	\ln \left(1-\left(2(1-\varepsilon)^{\tau / T}
	-  1\right)^{2}\right)}	
\end{equation*}
 is continuous on 
 $\mathbb{R}^{+}$
 with 
 $\lim _{\tau  \rightarrow 0^{+}} g(\tau)=\infty$, 
 and therefore there is a neighbourhood 
 of zero corresponding to 
 $\left(0, \tau (\varepsilon)\right)$ 
 within which 
 $ g $ is positive.

\section{Numerical simulations}
\label{sec:numerics}
We provide some numerical experiments to support the strong convergence proved in Theorem \ref{thm:convergence-rate-without-jump}. 
For convenience we recall here the original  Ait-sahalia-type rate model given in \eqref{eq:Ait-Sahalia-model}
\begin{equation*}\label{eq:Ait-Sahalia-model-Numerical}
	\dd X_t= (a_{-1} X_t^{-1} - a_{0}+ a_{1}
	X_t - a_{2} X_t^{\gamma} ) \, \dd t
	+ b X_t^{\theta} \dd W_t ,
	\quad
	t \in(0, T] ,
	\quad X_0 > 0
\end{equation*}
and the transformed SDE given in \eqref{eq:lamperti-translation1-jump}
\begin{equation*}
\begin{split}
d Y_t &	= f(Y_t) \dd t + b(1-\theta) \dd W_t,  \quad  t \in(0, T],\\
f(x) &  := (\theta-1) 
 \big( a_2 x^{\frac{\gamma-\theta}{1-\theta}}
 - a_1 x - 
a_{-1} x^{\frac{\theta+1} {\theta-1}} 
+  a_0 x^{\frac{\theta}{\theta-1}}
+ \frac{1}{2} b^2 \theta  x^{-1} \big).
\end{split}
\end{equation*}
We are interested in examining the performance of the splitting method \eqref{eq:splitting} (denoted \Splitting) and the tamed-splitting method \eqref{eq:tame-numerical-solution-without-jump} (denoted \TSM).
To construct a reference solution for our numerical experiments we apply a tamed Milstein method \cite{} to \eqref{eq:Ait-Sahalia-model} with a small step size $\tau=2^{-18}$.
We compare our methods to a backward Euler-Maruyama method applied to \eqref{eq:Ait-Sahalia-model}, denoted  \RefBEM. This was shown to have rate one half \cite{zhao2020backward}.
We also compare to other approximations of \eqref{eq:lamperti-translation1-jump} : a tamed Euler Maruyama method \cite{hutzenthaler2012} (denoted \TEM) and backward Euler Maruyama method (denoted \BEM).
In our experiments we take timesteps 
$\tau \in \{ 2^{-11},2^{-10},2^{-9},2^{-8},
	2^{-7}\}$ and fix $X_0=1$.
	For our convergence plots we take $1000$ realizations. In none of these realizations did \TSM give a negative value, that is at no point did we need to make use the backstop method \BEM in \eqref{eq:BEM_method}.
	Indeed the same was true when taking $10000$ realizations - illustrating the positivity Theorem \ref{thm:probability-positive}. 
	
\begin{figure}[tp]
	\includegraphics[width=0.5\linewidth, height=0.3\textheight]{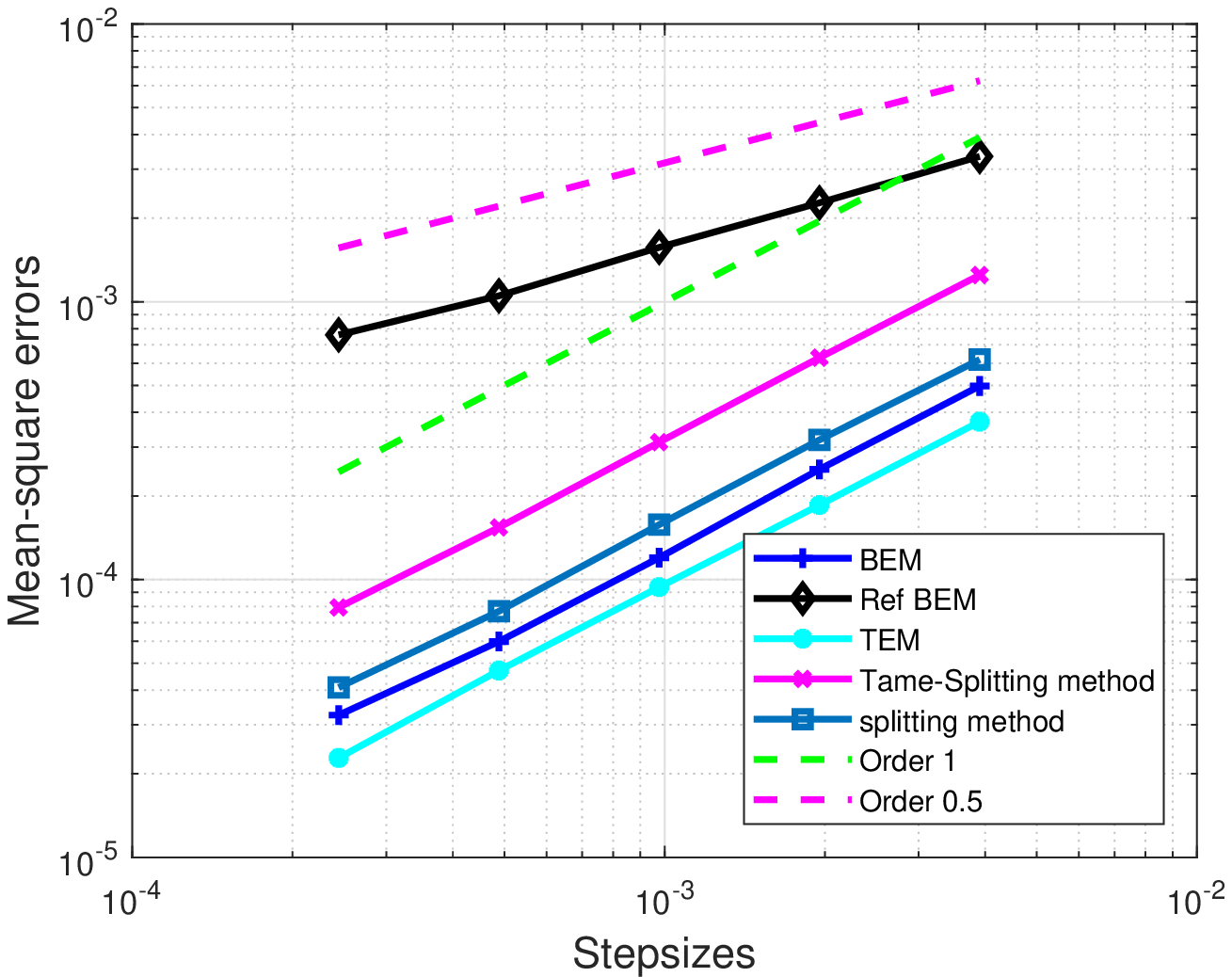}
	\includegraphics[width=0.5\linewidth, height=0.3\textheight]{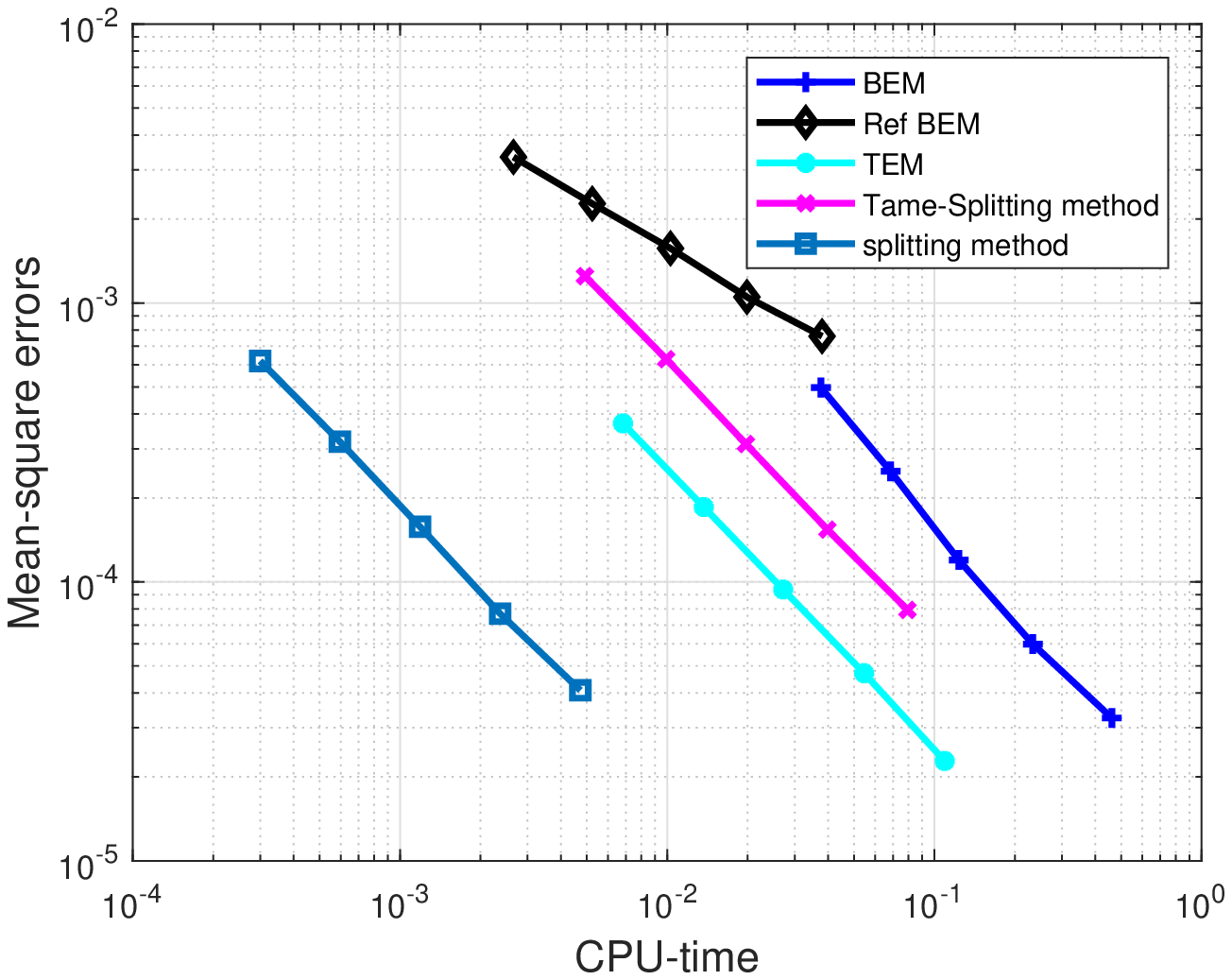}
	\caption[fig 1.]{Non-critical case. Convergence rates (left) and efficiency (right) of the numerical methods.}
	\label{fig1}
\end{figure}

\begin{figure}
	\includegraphics[width=0.5\linewidth, height=0.3\textheight]{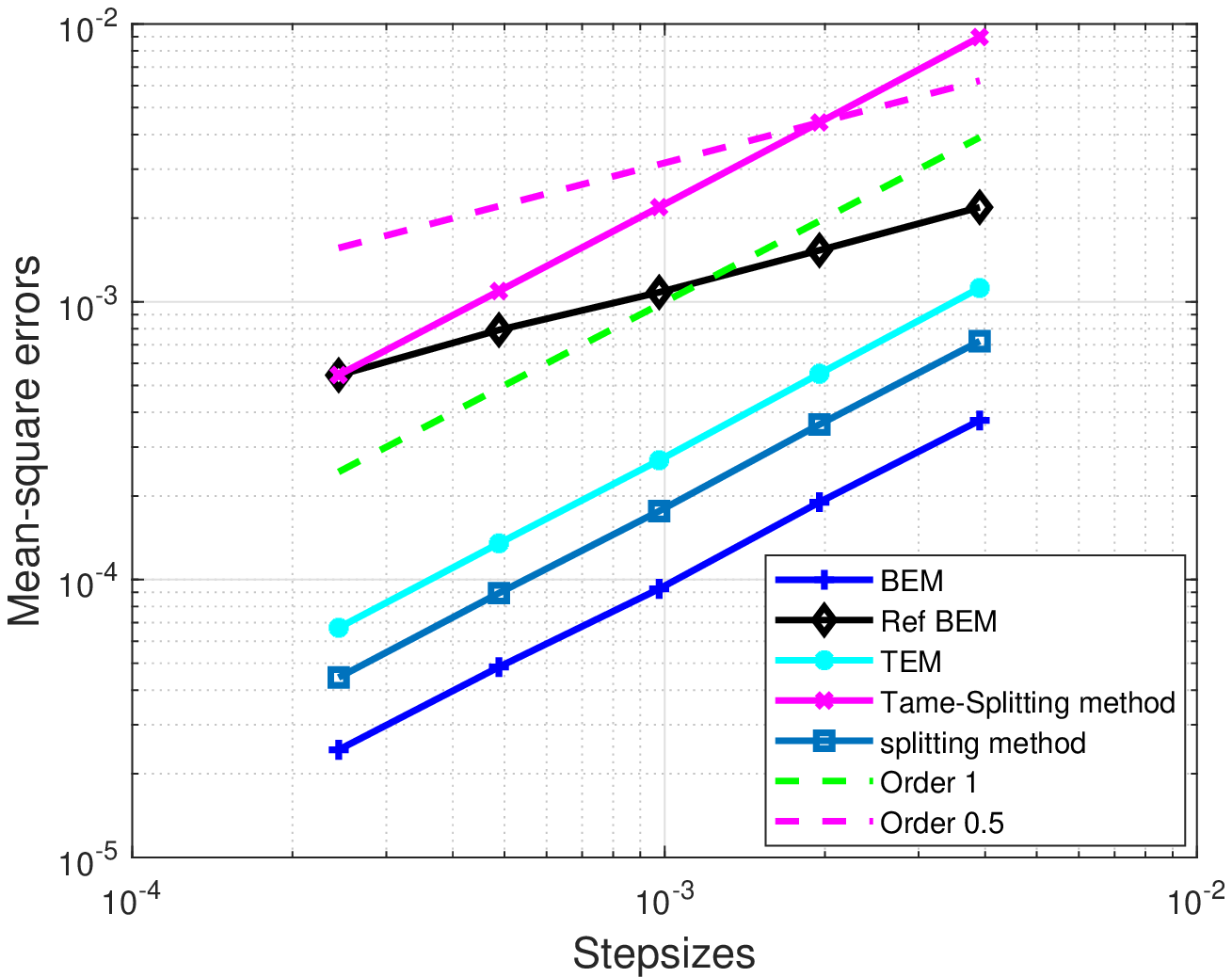}
	\includegraphics[width=0.5\linewidth, height=0.3\textheight]{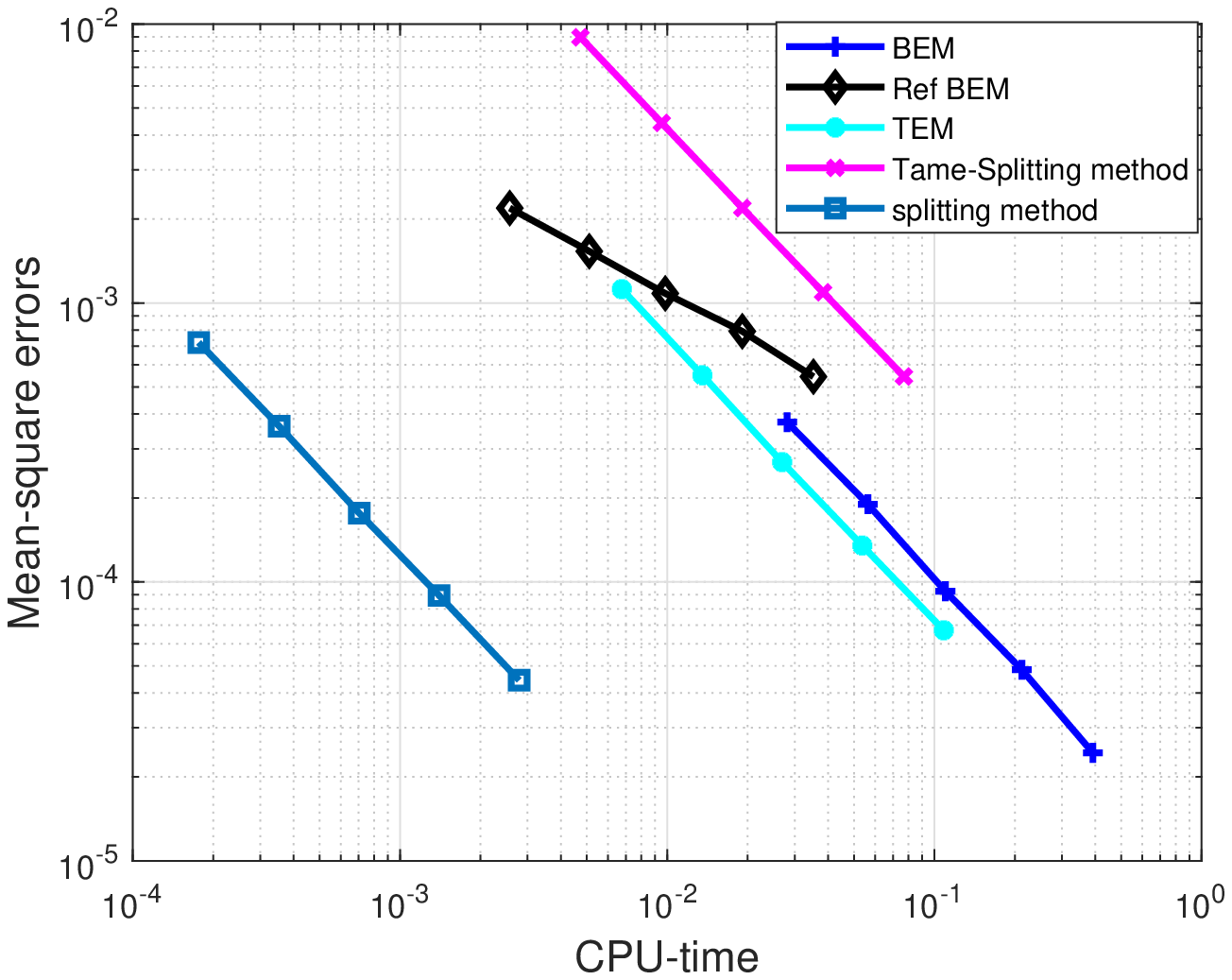}
	\caption[fig 2.]{Critical case. Convergence rates (left) and efficiency (right) of the numerical methods.}
	\label{fig2}
\end{figure}

We present simulations based on two sample sets that were originally derived in \cite{ait1996testing} from financial data. We fix 
	$a_{-1}=0.00107,a_0=0.0517,a_1=0.877,a_2=4.604,
	\gamma=3$ and $b=1$ and then take 
\begin{itemize}
	\item Non-critical case ( $ \gamma + 1 > 2 \theta $ ) : 
	$\theta=1.5. $
	\item Critical case ( $ \gamma + 1 = 2 \theta $ ) :
	$\theta=2.$
\end{itemize}
In Figure \ref{fig1} (left) we illustrate the convergence of our methods and the rates are given in Table \ref{tab:convergence}. We see that, as expected, \RefBEM has rate one half  and the other methods have rate 1. In this non-critical case we see that \TEM in fact has the best error constant. However, if we examine the efficiency (see Figure \ref{fig1} (right)) we see that \Splitting is the most efficient of the methods (and the two backward Euler Maruyama methods are the least efficient).
For the critical case we draw similar conclusions. 
In Figure \ref{fig2} (left) and Table \ref{tab:convergence} we illustrate convergence in this case. In this case it is surprising that the error constant for \TSM seems to increase by an order of magnitude. Examining the relative efficiency we again see that \Splitting is the most efficient by approximately two order of magnitude.
We illustrate the probability of numerical solutions take negative value.
It turns out that the probability of the numerical solutions being negative value is 0.
This means that our numerical method (\TSM) is almost positive. 
There is no need for backstop methods in numerical simulation.
Our numerical method (\TSM) is only possible to take negative values in theoretical analysis.  

The tamed-splitting method (\TSM) we analyse in this paper performs well in terms of the error constant but is slightly less efficient. However it is simpler to implement and can handle more complex situations as for general drift terms $F$ we may not exact solution of each equation by splitting.
%
%
%
\begin{table}[h]
\begin{center}
	\begin{tabular}{*{6}{c}}
		\toprule
		Parameter Sets  & \quad  \RefBEM & \BEM & 
		\quad  \TEM & \quad
		\TSM & \quad
		\Splitting \\
		\midrule
		Non-critical case &
		0.4909 & 0.9855 & 0.9795 & 0.9928
		& 0.9805  \\
		Critical case &
		0.5084  & 0.9880 & 1.0036 & 1.0084
		& 0.9852 \\
		\bottomrule
	\end{tabular}
	\caption{Rates of convergence estimated for the different methods, see Figs. \ref{fig1} and \ref{fig2}.}
	\label{tab:convergence}
	\end{center}
\end{table}

\bibliographystyle{siam}
\bibliography{Ait-Sahalia-V2}

\begin{thebibliography}{10}

\bibitem{ait1996testing}
{\sc Y.~Ait-Sahalia}, {\em Testing continuous-time models of the spot interest
  rate}, The review of financial studies, 9 (1996), pp.~385--426.

\bibitem{alfonsi2005discretization}
{\sc A.~Alfonsi}, {\em On the discretization schemes for the {CIR} (and
  {Bessel} squared) processes.}, Monte Carlo Methods \& Applications, 11
  (2005).

\bibitem{berkaoui2008euler}
{\sc A.~Berkaoui, M.~Bossy, and A.~Diop}, {\em Euler scheme for {SDEs} with
  non-{Lipschitz} diffusion coefficient: strong convergence}, ESAIM:
  Probability and Statistics, 12 (2008), pp.~1--11.

\bibitem{bossy2007efficient}
{\sc M.~Bossy and A.~Diop}, {\em An efficient discretisation scheme for one
  dimensional {SDEs} with a diffusion coefficient function of the form}, PhD
  thesis, INRIA, 2007.

\bibitem{conley1997short}
{\sc T.~G. Conley, L.~P. Hansen, E.~G. Luttmer, and J.~A. Scheinkman}, {\em
  Short-term interest rates as subordinated diffusions}, The Review of
  Financial Studies, 10 (1997), pp.~525--577.

\bibitem{Deng2018Generalized}
{\sc S.~{Deng}, C.~{Fei}, W.~{Fei}, and X.~{Mao}}, {\em Generalized
  {Ait-Sahalia-type} interest rate model with {Poisson} jumps and convergence
  of the numerical approximation}, Physica A., 533 (2019), p.~122057.

\bibitem{emmanuel2021truncated}
{\sc C.~Emmanuel and X.~Mao}, {\em Truncated {EM} numerical method for
  generalised {Ait}-{Sahalia}-type interest rate model with delay}, Journal of
  Computational and Applied Mathematics, 383 (2021), p.~113137.

\bibitem{gallant1997estimation}
{\sc A.~R. Gallant and G.~Tauchen}, {\em Estimation of continuous-time models
  for stock returns and interest rates}, Macroeconomic Dynamics, 1 (1997),
  pp.~135--168.

\bibitem{Gardon2004approximations}
{\sc A.~Gardo\'{n}}, {\em The order of approximations for solutions of
  {It\^{o}-type} stochastic differential equations with jumps}, Stoch Anal.
  Appl., 22 (2004), pp.~679--699.

\bibitem{giles2008multilevel}
{\sc M.~B. Giles}, {\em Multilevel {Monte} {Carlo} path simulation}, Operations
  research, 56 (2008), pp.~607--617.

\bibitem{higham2005convergence}
{\sc D.~J. Higham and X.~Mao}, {\em Convergence of {Monte Carlo} simulations
  involving the mean-reverting square root process}, Journal of Computational
  Finance, 8 (2005), pp.~35--61.

\bibitem{hong2005nonparametric}
{\sc Y.~Hong and H.~Li}, {\em Nonparametric specification testing for
  continuous-time models with applications to term structure of interest
  rates}, The Review of Financial Studies, 18 (2005), pp.~37--84.

\bibitem{hutzenthaler2012}
{\sc M.~Hutzenthaler, A.~Jentzen, and P.~E. Kloeden}, {\em Strong convergence
  of an explicit numerical method for {SDE}s with nonglobally {Lipschitz}
  continuous coefficients}, Ann. Appl. Probab., 22 (2012), pp.~1611--1641.

\bibitem{KLM2020}
{\sc C.~Kelly, G.~Lord, and H.~Maulana}, {\em The role of adaptivity in a
  numerical method for the {C}ox-{I}ngersoll- {R}oss model}, Journal of
  Computational and Applied Mathematics, 410 (2022), p.~114208.

\bibitem{kloeden_numerical_2011}
{\sc P.~Kloeden and E.~Platen}, {\em Numerical {Solution} of {Stochastic}
  {Differential} {Equations}}, Stochastic {Modelling} and {Applied}
  {Probability}, Springer Berlin Heidelberg, 2011.

\bibitem{kumar2019milstein}
{\sc C.~Kumar and S.~Sabanis}, {\em On milstein approximations with varying
  coefficients: the case of super-linear diffusion coefficients}, BIT Numerical
  Mathematics, 59 (2019), pp.~929--968.

\bibitem{lei2021first}
{\sc Z.~Lei, S.~Gan, and J.~Liu}, {\em First order strong approximation of
  {Ait}-{Sahalia}-type interest rate model with {Poisson} jumps}, arXiv
  preprint arXiv:2110.15482,  (2021).

\bibitem{lord2014introduction}
{\sc G.~J. Lord, C.~E. Powell, and T.~Shardlow}, {\em An introduction to
  computational stochastic PDEs}, vol.~50, Cambridge University Press, 2014.

\bibitem{Mao2008Stochastic}
{\sc X.~Mao}, {\em Stochastic Differential Equations and Applications},
  Horwood, 2008.

\bibitem{sasvari1999tight}
{\sc Z.~Sasvari and H.~Chen}, {\em Tight bounds for the normal distribution:
  10611}, The American Mathematical Monthly, 106 (1999), pp.~76--76.

\bibitem{Szpruch2011Numerical}
{\sc L.~Szpruch, X.~Mao, D.~J. Higham, and J.~Pan}, {\em Numerical simulation
  of a strongly nonlinear {Ait}-{Sahalia}-type interest rate model}, BIT
  Numerical Mathematics, 51 (2011), pp.~405--425.

\bibitem{zhang2017numerical}
{\sc Z.~Zhang and G.~E. Karniadakis}, {\em Numerical methods for stochastic
  partial differential equations with white noise}, Springer, 2017.

\bibitem{zhao2020backward}
{\sc Y.~Zhao, X.~Wang, and M.~Wang}, {\em On the backward {Euler} method for a
  generalized {Ait}-{Sahalia}-type rate model with {Poisson} jumps}, Numerical
  Algorithms,  (2020), pp.~1--21.

\end{thebibliography}
	
\end{document}